\documentclass[12pt]{amsart}  
 
\usepackage{amssymb} 
\usepackage{enumerate, amsfonts, latexsym,psfig,epsfig,amscd}

\newtheorem{theorem}{Theorem}[section]

\newtheorem{proposition}[theorem]{Proposition}

\newtheorem{lemma} [theorem]{Lemma} 
\newtheorem{corollary} [theorem] {Corollary} 

\newtheorem{remark}[theorem]{Remark}

\newtheorem{terminology}[theorem]{Terminology}

\theoremstyle{definition} 
\newtheorem{definition}[theorem]{Definition} 
\newtheorem{example}[theorem]{Example} 
 
\theoremstyle{remark}

\numberwithin{equation}{section} 
 
\newfont{\msb}{msbm10 scaled 1200} 
\newfont{\euf}{eufm10 scaled 1200}

\def\R{\mathbb R}
\def\Z{\mathbb Z}

\def\P{\mathbb P}

\def\A {\mathcal A}
\def\B {\mathcal B}

\def\SK {\underrightarrow{\text{Ker}}\, }

\catcode`\@=11 
\def\serieslogo@{\relax} 
\def\@setcopyright{\relax} 
\catcode`\@=12

\begin{document}

\title[CAT$(0)$ limit groups, II.]{Limits of (certain) CAT$(0)$ groups, II: The Hopf property and the shortening argument.}

% author  information 
\author[Daniel Groves]{Daniel Groves}
\address{Daniel Groves\\
Department of Mathematics \\
California Institute of Technology \\
Pasadena, CA, 91125, USA} 
\email{groves@caltech.edu} 

\date{August 4, 2004}

\begin{abstract}
This is the second in a series of papers about torsion-free groups which act properly and cocompactly on a CAT$(0)$ metric space with isolated flats and relatively thin triangles.  Our approach is to adapt the methods of Sela and others for word-hyperbolic groups to this context of non-positive curvature.

The main result in this paper is that (under certain technical hypotheses) such a group as above is Hopfian.  This (mostly) answers a question of Sela. 
\end{abstract}

\maketitle

\section{Introduction}
This paper is the second in a series, of which the first was \cite{CWIF1}.  Our approach is to consider methods and results of Sela and others in the context of negative curvature (word-hyperbolic groups), and adapt them to the context of non-positive curvature.  In particular, we consider a torsion-free group $\Gamma$ which acts properly, cocompactly and isometrically on a CAT$(0)$ metric space $X$ which has isolated flats and relatively thin triangles (as defined by Hruska, \cite{Hruska}).  For want of better terminology, we call $X$ a {\em CWIF space} and $\Gamma$ a {\em CWIF group}.\footnote{For technical reasons which arose in \cite{CWIF1}, we assume that the stabiliser in $\Gamma$ of any maximal flat in $X$ is free abelian;  such an action or group is called {\em toral}.}

In \cite{CWIF1} we provided an analogue of a construction of Paulin \cite{Paulin3}, by extracting from a sequence of pairwise non-conjugate homomorphisms $\{ h_n : G \to \Gamma \}$ a limiting action of $G$ on an $\R$-tree $T$.  This construction proceeds via a $G$-action on an asymptotic cone of $X$, and is briefly recalled in Section \ref{Prelim} below. 

Recall that a group $G$ is {\em Hopfian} if any surjective endomorphism $\phi : G \to G$ is an automorphism.  Sela \cite{SelaHopf} proved that torsion-free hyperbolic groups are Hopfian.  In contrast, Wise \cite{Wise} constructed a non-Hopfian CAT$(0)$ group.  Sela \cite[QI.8(i)]{SelaProblems} asked whether a group which acts properly and cocompactly by isometries on a CAT$(0)$ space with isolated flats is necessarily Hopfian.  The main result in this paper is that this is the case (under the additional hypotheses that the group is toral):

\begin{theorem} \label{GammaHopfian}
Suppose that $\Gamma$ is a torsion-free toral CWIF group.  Then $\Gamma$ is Hopfian.
\end{theorem}

Our approach to proving Theorem \ref{GammaHopfian} is to follow Sela's proof from \cite{SelaHopf} that torsion-free word-hyperbolic groups are Hopfian.  Other than the above-mentioned construction of an $\R$-tree, the main technical tool is the {\em shortening argument}.  This was developed in the context of automorphisms of word-hyperbolic groups by Rips and Sela in \cite{RipsSela1} and in the context of acylindrical accessibility by Sela in \cite{SelaAcyl}.  In the context of torsion-free toral CWIF groups, the shortening argument does not work in full generality; see Theorem \ref{NotShorteningArgument} below.  However, we are able to prove a strong enough version of it, Theorem \ref{ShorteningArgument}, in order to prove Theorem \ref{GammaHopfian}.

Given a finitely generated group $G$, the subgroup $\text{Mod}(G)$ (defined in Definition \ref{Mod} below) is generated by automorphisms which naturally arise from abelian splittings of $G$.  An immediate application of the shortening argument is the following

\begin{theorem} \label{ModinAut}
Suppose that $\Gamma$ is a freely indecomposable torsion-free toral CWIF group.  Then $\text{Mod}(\Gamma)$ has finite index in $\text{Aut}(\Gamma)$.
\end{theorem}

In the future work \cite{RelHypCSA} we will use the work of Drutu and Sapir \cite{DS} on asymptotic cones of relatively hyperbolic groups to extend the work of \cite{CWIF1} and the current paper to the context of torsion-free groups which are hyperbolic relative to a collection of free abelian groups.

The outline of this paper is as follows.  In Section \ref{Prelim} we recall the construction from \cite{CWIF1}.  We also recall the theory of stable isometric actions on $\R$-trees from \cite{BF} and \cite{SelaAcyl}, and the theory of abelian JSJ decompositions from \cite{RipsSela} (as slightly adapted for certain f.g. groups in \cite{Sela1}).  In Section \ref{NotShort} we describe the shortening argument and prove that it does not work for an arbitrary sequence of homomorphisms from a torsion-free toral CWIF group to itself (see Theorem \ref{NotShorteningArgument}).  In Section \ref{IET} we outline a version of the shortening argument that does work in this context and make some remarks, based on \cite{RipsSela1}.  We prove Theorem \ref{ShorteningArgument}, the shortening argument, assuming the technical results proved in Sections \ref{AxialSection} and \ref{Discrete}, and we also prove Theorem \ref{ModinAut}.  Sections \ref{AxialSection} and \ref{Discrete} are technical and contain the technical results needed to prove Theorem \ref{ShorteningArgument}.  Finally in Section \ref{HopfSection} we prove Theorem \ref{GammaHopfian}, by following and adapting the proof from \cite{SelaHopf}.

\section{Preliminaries} \label{Prelim}

In this section we gather the background results needed for this paper.  These include the construction from \cite{CWIF1}, and also the theory of isometric actions on $\R$-trees and that of JSJ decompositions.

For the definition of CAT$(0)$ spaces with isolated flats and/or relatively thin triangles, and for a host of motivating examples, see \cite{Hruska}.  We call such a CAT$(0)$ space a {\em CWIF space}, and a group which acts properly, cocompactly and isometrically on a CWIF space is called a {\em CWIF group}.   Recall that triangles in a CWIF space $X$ are relatively $\delta$-thin for some fixed constant $\delta > 0$ and that the definition of isolated flats involves a function $\phi : \R_+ \to \R_+$.  We use $\delta$ and $\phi$ later, but only to quote results from \cite{CWIF1}, so we do not repeat the appropriate definitions here.

\begin{remark}
Hruska and Kleiner have proved \cite{HruskaEmail} that the hypothesis that a cocompact CAT$(0)$ space with isolated flats has relatively thin triangles.  However, I do not know the proof, so we make the assumption that our space has both isolated flats and relatively thin triangles.
\end{remark}

For the remainder of this paper, fix the following notation.  The space $X$ is a CAT$(0)$ metric space which has isolated flats and relatively thin triangles.  The group $\Gamma$ is torsion-free and acts properly and cocompactly by isometries on $X$.  We also assume that the stabiliser in $\Gamma$ of a maximal flat in $X$ is free abelian (so the action is {\em toral}).  We also fix an arbitrary basepoint $x \in X$.

For much of the paper, we fix an arbitrary finitely generated group $G$, with finite generating set $\A$.

\subsection{From $\{ h_n : G \to \Gamma \}$, to $\mathcal C_\infty$, to the $\R$-tree $T$}

In this paragraph we recall the construction from \cite{CWIF1}.

\begin{terminology}
We say that two homomorphisms $h_1, h_2 : G \to \Gamma$ are {\em conjugate} if there is $\gamma \in \Gamma$ so that $h_1 = \tau_\gamma \circ h_2$, where $\tau_\gamma$ is the inner automorphism of $\Gamma$ induced by $\gamma$.  Otherwise, $h_1$ and $h_2$ are {\em non-conjugate}.
\end{terminology}

\begin{definition}
Suppose that $\{ h_n : G \to \Gamma \}$ is a sequence of homomorphisms.  The {\em stable kernel} of $\{ h_n \}$, denoted $\SK (h_n)$, is the set of all $g \in G$ so that $g \in \text{ker}(h_n)$ for all but finitely many $n$.
\end{definition}

From a sequence of pairwise non-conjugate homomorphisms $\{ h_n : G \to \Gamma \}$, it is straightforward to construct an isometric action of $G$ on $X_\omega$, where $\omega$ is a non-principal ultrafilter, and $X_\omega$ is an asymptotic cone of $X$.\footnote{For the purposes of this paper, we are unconcerned which ultrafilter $\omega$ is, but the sequence of scaling factors are important; see \cite{CWIF1} for the details of this construction in the context of this paper.}  The basepoint $x_\omega$ of $X_\omega$ is the point represented by the constant sequence $\{ x \}$.  A key feature of the action of $G$ on $X_\omega$ is that there is no global fixed point.  This construction is carried out in \cite{Paulin3} when $\Gamma$ is hyperbolic, and in \cite{KL} when $\Gamma$ is CAT$(0)$.  

When $X$ has isolated flats and relatively thin triangles, the properties of $X_\omega$ were enumerated in \cite{CWIF1}.  Since \cite{CWIF1} was written, the paper of Dru\c{t}u and Sapir \cite{DS} has appeared as a preprint.  In the terminology of \cite{DS}, $X_\omega$ is a {\em tree-graded} metric space.  Using the results of \cite{DS}, in the future work \cite{RelHypCSA} the results from this paper and from \cite{CWIF1} are generalised to the case where $\Gamma$ is a torsion-free group which is hyperbolic relative to a collection of free abelian subgroups.

Given the action of $G$ on $X_\omega$, and the basepoint $x_\omega \in X_\omega$, define the space $\mathcal C_\infty \subset X_\omega$ to be the union of the geodesics $[x_\omega , g . x_\omega]$, for $g \in G$, along with any flats $E \subset X_\omega$ which contains a non-degenerate open triangle contained in some triangle $\Delta (g_1 .  x_\omega, g_2 . x_\omega, g_3 . x_\omega)$, where $g_1,g_2,g_3 \in G$.  With the induced path metric, $\mathcal C_\infty$ is a tree-graded metric space, is convex in $X_\omega$, is $G$-invariant, and the $G$-action on $\mathcal C_\infty$ does not have a global fixed point.  Also, $\mathcal C_\infty$ is seperable.  It is the seperability which allows the proof of the following lemma.  Before stating the lemma, we recall the $G$-equivariant Gromov topology on metric spaces equipped with isometric $G$-actions.

Suppose that $\{ (Y_n, y_n, \lambda_n) \}_{n=1}^\infty$ and $(Y,y,\lambda)$ are triples consisting of metric spaces, basepoints, and actions $\lambda_n : G \to \text{Isom}(Y_n)$, $\lambda : G \to \text{Isom}(Y)$.  Then $(Y_n,y_n, \lambda_n) \to (Y,y,\lambda)$ in the {\em $G$-equivariant Gromov topology} if and only if:  for any finite subset $K$ of $Y$, any $\epsilon > 0$ and any finite subset $P$ of $G$, for sufficiently large $n$ there are subsets $K_n$ of $Y_n$ and bijections $\rho_n : K_n \cup \{ y_n \} \to K \cup \{ y \}$ such that $\rho(y_n) = y$ and for all $s_n, t_n \in K_n \cup \{ y_n \}$ and all $g_1, g_2 \in P$ we have
\[	|d_Y(\lambda (g_1) . \rho_n(s_n) , \lambda(g_2) . \rho_n(t_n)) - d_{Y_n}(\lambda_n (g_1) . s_n , \lambda_n(g_2), t_n ) | < \epsilon .	\]

Given a homomorphism $h : G \to \Gamma$, we define the {\em length} of $h$:
\[	\| h \| = \max_{g \in \A}d_X(x, h(g) . x) .	\]
We associate to $h$ a triple $(X_h, x_h, \lambda_h)$ as follows:  let $X_h$ be the convex hull in $X$ of $G . x$, endowed with the metric $\frac{1}{\| h \| }d_X$, let $x_h = x$, and let $\lambda_h = \iota \circ h$, where $\iota : \Gamma \to \text{Isom}(X)$ is the (fixed) homomorphism given by the action of $\Gamma$ on $X$.

\begin{lemma} \cite[Lemma 3.15]{CWIF1}
Let $\Gamma, X$ and $G$ be as described above, and let $\{h_n : G \to \Gamma \}$ be a sequence of pairwise non-conjugate homomorphisms.  Let $X_\omega$ be the asymptotic cone of $X$, let $x_\omega$ be the basepoint on $X_\omega$ and let $\mathcal C_\infty$ be as described above.  Let $\lambda_\infty : G \to \text{Isom}(\mathcal C_\infty)$ denote the action of $G$ on $\mathcal C_\infty$ and $(\mathcal C_\infty, x_\omega, \lambda_\infty)$ the associated triple.

Then there exists a subsequence $\{ f_i \} \subseteq \{ h_i \}$ so that the triples $(X_{f_i}, x_{f_i}, \lambda_{f_i})$ converge to $(\mathcal C_\infty, x_\omega, \lambda_\infty)$ in the $G$-equivariant Gromov topology.
\end{lemma}

In Section 4 of \cite{CWIF1}, we describe how to extract an $\R$-tree $T$ from $\mathcal C_\infty$ so that $T$ has an isometric $G$-action with no global fixed point and so that the kernel of the $G$-action on $T$ is exactly the kernel of the $G$-action on $\mathcal C_\infty$.  The idea is to take each flat $E$ in $\mathcal C_\infty$ and project it to a line $p_E$ (there are a number of conditions on which lines are allowed).  Since the action of $\text{Stab}(E)$ on a maximal flat $E \subset \mathcal C_\infty$ is by translations, this action projects to an action of $\text{Stab}(E)$ on $p_E$ by translations.  The collection of the lines $p_E$ is denoted by $\P$.  There is a natural map, which we call {\em projection}, from $\mathcal C_\infty \to T$ which restricts to a bijection from $\mathcal C_\infty \smallsetminus \{ E \mid \mbox{$E$ is a maximal flat in $\mathcal C_\infty$} \}$ to $T \smallsetminus \P$.  See Section 4 of \cite{CWIF1} for the details of the construction of $T$.

The main construction of \cite{CWIF1} is summarised in the following

\begin{theorem} [cf. Theorem 4.4, \cite{CWIF1}] \label{Texists}
Suppose that $\Gamma$ is a torsion-free toral CWIF group and that $G$ is a finitely generated group.  Let $\{ h_n : G \to \Gamma \}$ be a sequence of pairwise non-conjugate homomorphisms.  There is a subsequence $\{ f_i \}$ of $\{ h_i \}$ and an action of $G$ on an $\R$-tree $T$ so that if $K$ is the kernel of the action of $G$ on $T$ and $L := G/K$ then
\begin{enumerate}
\item The stabiliser in $L$ of any nondegenerate segment in $T$ is free abelian;
\item If $T$ is isometric to a real line then $L$ is free abelian, and for all but finitely $n$ the group  $h_n(G)$ is free abelian;
\item If $g \in G$ stabilises a tripod in $T$ then $g \in \SK (f_i) \subseteq K$;
\item Let $[y_1, y_2] \subset [y_3, y_4]$ be non-degenerate segments in $T$, and assume that $\text{Stab}_L([y_3,y_4])$ is nontrivial.  Then
\[	\text{Stab}_L([y_1,y_2]) = \text{Stab}_L([y_3,y_4]).	\]
In particular, the action of $L$ on $T$ is stable; and
\item $L$ is torsion-free.
\end{enumerate}
\end{theorem}

\begin{definition}
Suppose that $\Gamma$ is a torsion-free toral CWIF group.  A {\em $\Gamma$-limit group} is a finitely generated group $L$ which is either isomorphic to a subgroup of $\Gamma$ or is of the form $L \cong G/K$, where $G$ is a finitely generated group, and $K$ is the kernel of the action of $G$ on $\mathcal C_\infty$, where $\mathcal C_\infty$ arises as above from a sequence of pairwise non-conjugate homomorphisms $\{ h_n : G \to \Gamma \}$.  Those groups of the second kind are called {\em strict $\Gamma$-limit groups}.
\end{definition}

\begin{definition}
The group $H$ is {\em commutative transitive} if for all elements $u_1, u_2, u_3 \in H \smallsetminus \{ 1 \}$, whenever $[u_1,u_2] = 1$ and $[u_2, u_3] = 1$, necessarily $[u_1, u_3] = 1$.

A subgroup $K$ of a group $H$ is said to be {\em malnormal} if for all $h \in H \smallsetminus K$ we have $hKh^{-1} \cap K = \{ 1 \}$.  The group $H$ is {\em CSA} if any maximal abelian subgroup of $H$ is malnormal.
\end{definition}

\begin{remark}
By \cite[Lemma 2.19]{CWIF1}, a torsion-free CWIF group is toral if and only if it is CSA.  Since the phrase `CSA CWIF group' is even worse than the terminology in this paper, we persist with the term `toral'.
\end{remark}

\begin{lemma} \cite[Corollary 5.7]{CWIF1} \label{LimitCSA}
Suppose that $\Gamma$ is a torsion-free toral CWIF group and that $L$ is a $\Gamma$-limit group.  Then $L$ is commutative transitive and CSA.
\end{lemma}

\subsection{Isometric actions on $\R$-trees}

In this paragraph we recall a result of Sela from \cite{SelaAcyl}.  Given a finitely generated group $G$ and an $\R$-tree $T$ with an isometric $G$-action, Theorem \ref{RtreeStructure} below gives a decomposition of $T$ which induces a graph of groups decomposition of $G$.  In the case that $G$ is finitely presented, this result follows immediately from Rips Theory; see Bestvina and Feighn, \cite{BF}.

There are two sets of terminology for the components of the above-mentioned decomposition.  Since we are quoting Sela's result, we use his (Rips') terminology.  However, we assume that all axial components are isometric to a real line.  Using Rips and Sela's definition of axial (see \cite[\S 10]{RipsSela1}), one other case could arise in the arguments that follow (where our group splits as $A \ast_{[a,b]}\langle a,b \rangle$).  Just as noted in \cite[\S4, p.346]{RipsSela1}, we can treat this case as an IET component.  Thus, without further mention, we consider all axial components to be isometric to a real line.

The following theorem of Sela is used to decompose our limiting $\R$-trees.

\begin{theorem} [Theorem 3.1, \cite{SelaAcyl}; see also Theorem 1.5, \cite{Sela1}] \label{RtreeStructure}
Let $G$ be a freely indecomposable finitely generated group which admits a stable isometric action on a real tree $Y$.  Assume that the stabiliser in $G$ of each tripod in $T$ is trivial.
\begin{enumerate}
\item[1)] There exist canonical orbits of subtrees of $T$, denoted $T_1, \ldots , T_k$, with the following properties:
\begin{enumerate}
\item[(i)] for each $g \in G$ and each $i, j \in \{ 1, \ldots , k \}$ with $i \neq j$, the subtree $g.T_i$ intersects $T_j$ in at most a single point;
\item[(ii)] for each $g \in G$ and each $i \in \{ 1, \ldots , k \}$, the subtree $g.T_i$ is either equal to $T_i$ or intersects $T_i$ in at most a single point;
\item[(iii)] The action of $\text{Stab}_G(T_i)$ on $T_i$ is of axial or IET type;
\end{enumerate}
\item[2)] $G$ is the fundamental group of a graph of groups with:
\begin{enumerate}
\item[(i)] Vertices corresponding to orbits of branching points with non-trivial stabliser in $T$;
\item[(ii)] Vertices corresponding to the orbits of the canonical subtrees $T_1, \ldots , T_k$ which are of axial or IET type.  The groups associated with these vertices are conjugates of the stabilisers of these subtrees.  To a stabiliser of an IET component there exists an associated $2$-orbifold, $\mathcal O$.  Any element of $\pi_1(\mathcal O)$ which corresponds to a boundary component or branching point in $\mathcal O$ stabilises a point in $T$.  For each stabiliser of an IET subtree we add edges that connect the vertex stabilised by it and the vertices stabilised by its boundary components and branching points;
\item[(ii)] Edges corresponding to orbits of edges between branching points with non-trivial stabiliser in the discrete part of $T$ (see Terminology \ref{DiscretePart} below) with edge groups which are conjugates of the stabilisers of these edges;
\item[(iii)] Edges corresponding to orbits of points of intersection between the orbits of $T_1, \ldots , T_k$.
\end{enumerate}
\end{enumerate}
\end{theorem}

\begin{terminology} \label{DiscretePart}
Let $G$ and $T$ be as in Theorem \ref{RtreeStructure} above.  The {\em discrete part of $T$} is the union of the metric closures of the connected components of $T \setminus \left( \bigcup\limits_{i=1}^k G T_i \right)$.
\end{terminology}

In order to prove Theorem \ref{GammaHopfian}, we also need the following 

\begin{remark} \cite[Remark 3.15, p.15]{BFSela} \label{RelativeRtreeStructure}
Theorem \ref{RtreeStructure} holds more generally if the assumption that $G$ is freely indecomposable is replaced by the assumption that $G$ is freely indecomposable rel point stabilisers; i.e. if $\mathcal V$ is the subset of $G$ of elements acting elliptically on $T$, then $G$ cannot be expressed non-trivially as $A \ast B$ with all $g \in \mathcal V$ conjugate into $A \cup B$.
\end{remark}

\begin{remark} \label{NoThin}
In the theory of stable isometric actions on $\R$-trees, there is one further type of component arising in the decomposition of $T$.  This is called `thin' in \cite{BF} and was discovered and investigated by Levitt (see \cite{Levitt}).  However, in case $G$ is freely indecomposable and the stabiliser of any non-degenerate tripod is trivial (both of these conditions hold in all cases in this paper), thin components do not arise.
\end{remark}

\subsection{Acylindrical accessibility} \label{AcylSection}

The theory of JSJ decompositions is used in the next paragraph to refine the graph of groups decomposition given in Theorem \ref{RtreeStructure} to a `canonical' graph of groups decomposition.  For finitely presented groups, we can use the accessibility result from \cite{BF2}.   However, the groups we consider need not be finitely presented, so we need to use the {\em acylindrical accessibility} of \cite{SelaAcyl}.

\begin{definition} \cite{SelaAcyl}
A splitting of a group $H$ is {\em reduced} if the label of every vertex of valence $2$ properly contains the labels of both edges incident to it.  Let $Y$ be the Bass-Serre tree for a given splitting of the group $H$.  We say that the splitting, and $Y$, are {\em $k$-acylindrical} if they are reduced, $T$ is minimal, and for all $h \in H \smallsetminus \{ 1 \}$ the fixed set of $h$ in $Y$ has diameter at most $k$.
\end{definition}

\begin{theorem} [\cite{SelaAcyl}; see also \cite{Weidmann}]
Let $H$ be a finitely generated, freely indecomposable group.  For a given $k$ there exists an integer $\xi(k,H)$ so that the number of vertices and edges in any $k$-acylindrical splitting of $H$ does not exceed $\xi(k,H)$.
\end{theorem}

\begin{lemma} [cf. Lemma 2.1, p.39, \cite{Sela1}] \label{AbelianElliptic}
Let $\Gamma$ be a torsion-free toral CWIF group and $L$ a $\Gamma$-limit group.  Suppose that $M$ is a noncyclic maximal abelian subgroup of $L$, and $A$ is an abelian subgroup of $L$.  Then:
\begin{enumerate}
\item If $L = U \ast_A V$ then $M$ can be conjugated into either $U$ or $V$;
\item If $L = U \ast_A$ then either $M$ can be conjugated into $U$, or $L = U \ast_A M'$, where $M'$ is a conjugate of $M$.
\end{enumerate}
\end{lemma}
\begin{proof}
Given Lemma \ref{LimitCSA}, the proof is identical to that of \cite[Lemma 2.1, p.39]{Sela1}.
\end{proof}

Exactly as in \cite[\S2, p.40]{Sela1}, we have the following

\begin{lemma} [cf. Lemma 2.3, p.40, \cite{Sela1}]
Let $\Gamma$ be a torsion-free toral CWIF group and $L$ a $\Gamma$-limit group.  A splitting of $L$ in which all edge groups are abelian and all noncyclic abelian groups are elliptic can always be modified (by modifying boundary monomorphisms by conjugations and sliding operations) to be $2$-acylindrical.
\end{lemma}

Using the above two lemmas, we can ensure that all of our splittings are $2$-acylindrical.

\subsection{JSJ decompositions}

There are a number of approaches to JSJ decompositions -- see \cite{RipsSela}, \cite{DunSag}, \cite{FP} and \cite{SS}.  We follow the approach of Rips and Sela, adapted as in \cite{SelaHopf} and \cite{Sela1} to the context of certain groups which are not necessarily finitely presented.

In the following theorem, the phrase `under consideration' refers to $2$-acylindrical abelian splittings where all non-cyclic abelian subgroups are elliptic.

\begin{theorem} [see Theorem 2.7, \cite{Sela1}] \label{AbelianJSJ}
Suppose that $L$ is a freely indecomposable strict $\Gamma$-limit group. There exists a reduced unfolded splitting of $L$ with abelian edge groups, which we call an {\em abelian JSJ decomposition of $L$} with the following properties:
\begin{enumerate}
\item Every canonical maximal QH (CMQ) subgroup of $L$ is conjugate to a vertex group in the JSJ decomposition.  Every QH subgroup of $L$ can be conjugated into one of the CMQ subgroups of $L$.  Every vertex group in the JSJ decomposition of $L$ which is not a CMQ subgroup is elliptic in any abelian splitting of $L$ under consideration; \label{JSJ-Vertex}
\item A one edge abelian splitting under consideration, $L = D \ast_A E$ or $L = D \ast_A$, which is hyperbolic in another elementary abelian splitting is obtained from the abelian JSJ decomposition of $L$ by cutting a $2$-orbifold corresponding to a CMQ subgroup of $L$ along a weakly essential s.c.c;
\item Let $\Upsilon$ be a one edge abelian splitting under consideration -- $L = D \ast_A E$ or $L = D \ast_A$ -- and suppose that $\Upsilon$ is elliptic with respect to any other one edge abelian splitting of $L$ under consideration.  Then $\Upsilon$ is obtained from the JSJ decomposition of $L$ by a sequence of collapsings, foldings and conjugations;
\item If JSJ$_1$ is another JSJ decomposition of $L$, then JSJ$_1$ is obtained from the JSJ decomposition by a sequence of slidings, conjugations and modifying boundary monomorphisms by conjugations.
\end{enumerate}
\end{theorem}

See \cite{RipsSela} for details on the undefined terms in the above theorem.

\section{Not the shortening argument} \label{NotShort}

\begin{definition} [Dehn twists]
Let $G$ be a finitely generated group.  A {\em Dehn twist} is an automorphism of one of the following two types:
\begin{enumerate}
\item Suppose that $G = A \ast_C B$ and that $c$ is contained in the centre of $C$.  Then define $\phi \in \text{Aut}(G)$ by $\phi(a) = a$ for $a \in A$ and $\phi(b) = cbc^{-1}$ for $b \in B$;
\item Suppose that $G = A \ast_C$, that $c$ is in the centre of $C$, and that $t$ is the stable letter of this HNN extension.  Then define $\phi \in \text{Aut}(G)$ by $\phi(a) = a$ for $a \in A$ and $\phi(t) = tc$.
\end{enumerate}
\end{definition}

\begin{definition} [Generalised Dehn twists]
Suppose $G$ has a graph of groups decomposition with abelian edge groups, and $A$ is an abelian vertex group in this decomposition.  Let $A_1 \le A$ be the subgroup generated by all edge groups connecting $A$ to other vertex groups in the decomposition.  Any automorphism of $A$ that fixes $A_1$ elementwise can be naturally extended to an automorphism of the ambient group $G$.  Such an automorphism is called a {\em generalised Dehn twist} of $G$.
\end{definition}

\begin{definition} \label{Mod}
Let $G$ be a finitely generated group.  We define $\text{Mod}(G)$ to be the subgroup of $\text{Aut}(G)$ generated by:
\begin{enumerate}
\item Inner automorphisms;
\item Dehn twists arising from splittings of $G$ with abelian edge groups; and
\item Generalised Dehn twists arising from graph of groups decompositions of $G$ with abelian edge groups.
\end{enumerate}
\end{definition}
Similar definitions are made in \cite[\S 5]{Sela1} and \cite[\S 1]{BFSela}

Suppose that $\Gamma$ is a torsion-free toral CWIF group with a toral, proper and cocompact action by isometries on a CWIF space $X$, with basepoint $x \in X$.  Suppose also that $G$ is a finitely generated group, with finite generating set $\A$.  Let $h : G \to \Gamma$ be a homomorphism.  Recall that in Section \ref{Prelim} we defined the {\em length} of $h$ by
\[	\| h \| := \max_{g \in \A} \big\{ d_X(x, h(g) . x) \big\} . 	\]

\begin{definition} [cf. Definition 4.2, \cite{BFSela}] \label{ShortHomo}
We define an equivalence relation on the set of homomorphisms $h : G \to \Gamma$ by setting $h_1 \sim h_2$ if there is $\alpha \in \text{Mod}(G)$ and $\gamma \in \Gamma$ so that $h_1 = \tau_\gamma \circ h_2 \circ \alpha$, where $\tau_\gamma$ is the inner automorphism of $\Gamma$ induced by $\gamma$.

A homomorphism $h : G \to \Gamma$ is {\em short} if for any $h'$ such that $h \sim h'$ we have $\| h \| \le \| h' \|$.
\end{definition}

The following is one of the main technical results of this paper.

\begin{theorem} [Shortening Argument] \label{ShorteningArgument}
Suppose that $\Gamma$ is a non-abelian, freely indecomposable, torsion-free toral CWIF group, and suppose that the sequence of automorphisms $\{ h_n : \Gamma \to \Gamma \}$ converges to an action $\eta : \Gamma \to \text{Isom}(\mathcal C_\infty)$ as above.  Then for all but finitely many $n$ the homomorphism $h_n$ is not short.
\end{theorem}

One may hope to prove Theorem \ref{ShorteningArgument} for arbitrary sequences of homomorphisms, rather than just sequences of automorphisms, and where the homomorphisms are from an arbitrary finitely generated group $G$ to $\Gamma$, rather than from $\Gamma$ to itself.  However, this is impossible, as shown by the next result.  It is unclear exactly what conditions to impose on a sequence of homomorphisms from a finitely generated group $G$ to $\Gamma$ to allow the Shortening Argument to be applied.

\begin{theorem} \label{NotShorteningArgument}
There exists a torsion-free toral CWIF group $\Gamma$, and a sequence of short homomorphisms $h_n : \Gamma \to \Gamma$ which converge to a faithful action of $\Gamma$ on a limiting space $\mathcal C_\infty$.
\end{theorem}

Although we can produce explicit $\Gamma$, we do not construct a specific sequence of short homomorphisms.  Rather we argue by contradiction.  Theorem \ref{NotShorteningArgument} is immediately implied by Proposition \ref{FinConj} and Example \ref{NotcoHopf} below.

\begin{proposition} [cf. Theorem 7.1, p.364, \cite{RipsSela1}] \label{FinConj}
Suppose that there is a freely indecomposable, nonabelian and torsion-free toral CWIF group $\Gamma$ such that there are infinitely many conjugacy classes of embeddings of $\Gamma$ into itself.  Then Theorem \ref{NotShorteningArgument} holds for this choice of $\Gamma$.
\end{proposition}
\begin{proof}
We follow the proof of \cite[Theorem 7.1]{RipsSela1}.  

Suppose that there are infinitely many conjugacy classes of embeddings of $\Gamma$ in $\Gamma$.  Then there exists a sequence of embeddings $\{ \rho_i : \Gamma \to \Gamma \}$ so each that $\rho_i$ comes from a different conjugacy class and each $\rho_i$ is short.

By Theorem \ref{Texists}, there is a subsequence converging to an action of $\Gamma$ on an $\R$-tree $T$.  Since $\Gamma$ is not abelian and each $\rho_i$ is an embedding, $T$ is not isometric to a real line, so the kernel of the action of $\Gamma$ on $T$ is $\SK (\rho_i)$, which is the identity.  Therefore the action of $\Gamma$ on $T$ is faithful.  This implies that Theorem \ref{NotShorteningArgument} holds for this choice of $\Gamma$.
\end{proof}

\begin{example} \label{NotcoHopf}
Let $H$ be a freely indecomposable non-elementary and torsion-free CAT$(-1)$ group, and let $w \in H$ be an element such that the normaliser of $w$ is $\langle w \rangle$.  Let $\Gamma = H \ast_{\Z} = \langle H, t\ |\ tw t^{-1} = w \rangle$.  

It is easy to see that $\Gamma$ is a torsion-free toral CWIF group, which is also finitely presented, freely indecomposable and non-abelian.

However, $\Gamma$ contains infinitely many conjugacy classes of embeddings of $\Gamma$.  Define $\phi_n : \Gamma \to \Gamma$ by $\phi_n(h) = h$ for all $h \in H$ and $\phi_n(t) = t^n$.  It is easy to see that each $\phi_n$ is an embedding of $\Gamma$ into itself.  Since the relations of $\Gamma$ and conjugation both preserve the exponent sum of $t$, it is also clear that the $\phi_n(\Gamma)$ are all distinct, even up to conjugacy.
\end{example}

Note that in the above example, $\Gamma$ is not co-Hopfian\footnote{Recall that a group is {\em co-Hopfian} if every injective endomorphism is surjective.}, in contrast to the result of Sela \cite{SelaGAFA} that torsion-free non-elementary hyperbolic groups are co-Hopfian if and only if they are freely indecomposable.  

Example \ref{NotcoHopf} was suggested by examples of Kleinian groups in \cite{OP}.  In \cite{DP} the question of when a Kleinian group is not co-Hopfian was investigated.  

In the future work \cite{coHopf} we examine the question of which groups which are hyperbolic relative to free abelian groups are co-Hopfian (this class is more general than torsion-free toral CWIF groups).  In \cite{Dahmani}, Dahmani investigates the question of the finiteness of conjugacy classes of embeddings of subgroups of relatively hyperbolic groups (without the condition that the parabolic subgroups are abelian).

\section{The shortening argument -- introduction} \label{IET}

In this section we outline the proof of Theorem \ref{ShorteningArgument} (the complete proof is contained in this and the subsequent two sections):

\medskip

{\bf Theorem \ref{ShorteningArgument}} (Shortening Argument).
{\em Suppose that $\Gamma$ is a nonabelian, freely indecomposable, torsion-free toral CWIF group, and suppose that the sequence of automorphisms $\{ h_n : \Gamma \to \Gamma \}$ converges to an action $\eta : \Gamma \to \text{Isom}(\mathcal C_\infty)$ as above.  Then for all but finitely many $n$ the homomorphism $h_n$ is not short.
}

\medskip

\begin{remark}
Although we call the above theorem the `Shortening Argument', at least for hyperbolic groups it is more a collection of ideas which can be used in myriad situations.  The above theorem is enough to prove Theorem \ref{ModinAut}, but we need to refer to the proof rather than the statement of Theorem \ref{ShorteningArgument} in order to prove Theorem \ref{GammaHopfian}.  Our hope is that the shortening argument in some form will be applicable to many further problems about CWIF groups and relatively hyperbolic groups.
\end{remark}

Let $\{ h_n : \Gamma \to \Gamma \}$ be a sequence of pairwise non-conjugate automorphisms.  Since $\Gamma$ is non-abelian, the action of $\Gamma$ on the limiting space $\mathcal C_\infty$ is faithful, and the action of $\Gamma$ on the associated $\R$-tree $T$ is also faithful.  We prove that for all but finitely many $n$, the homomorphism $h_n$ is not short.

Since the action of $\Gamma$ on $T$ is faithful, $\Gamma$ is itself a strict $\Gamma$-limit group, and by Theorem \ref{Texists} the stabiliser in $\Gamma$ of any tripod in $T$ is trivial.

The approach to proving Theorem \ref{ShorteningArgument} is as follows:  we consider the finite generating set $\A_1$ of $\Gamma$, and the basepoint $y$ of $T$.  We consider the paths $[y, u.y]$ where $u \in \A_1$.  These paths can travel through various types of subtrees of $T$; the IET subtrees, the axial subtrees, and the discrete part of $T$.\footnote{Recall by Remark \ref{NoThin} that there are no thin components in the decomposition of $T$.}  Depending on the types of subtrees which have positive length intersection with $[y,u.y]$, we need various types of arguments which allow us to shorten the homomorphisms which `approximate' the action of $\Gamma$ on $T$.

Mostly, we follow the shortening argument as developed in \cite{RipsSela1}.  There are two main obstacles to implementing this strategy in the context of torsion-free toral CWIF groups.  Note that the automorphisms $h_n : \Gamma \to \Gamma$ actually `approximate' the action of $G$ on $\mathcal C_\infty$, from which the action of $\Gamma$ on $T$ is extracted.  The two main problems are: (i) those lines $p_E \in \P$ which correspond to flats $E \in \mathcal C_\infty$; and (ii) that triangles in the approximating spaces are only relatively thin, not actually thin.

\subsection{IET components}

The following theorem of Rips and Sela deals with IET components.

\begin{theorem} \cite[Theorem 5.1, pp. 346-347]{RipsSela1} \label{IETTheorem}
Let $G$ be a finitely presented freely indecomposable group and assume that $G \times T \to T$ is a small stable action of $G$ on a real tree $T$ with trivial stabilisers of tripods.  Let $U$ be a finite subset of $G$ and let $y \in T$.  Then there exists $\phi_I \in \text{Mod}(G)$ such that for any $u \in U$, if $[y,u(y)]$ has an intersection of positive length with some IET-component of $T$ then:
\[	d_Y(y,\phi_I(u) . y) <d_T(y,u(y)),	\]
and otherwise $\phi_I(u) = u$.
\end{theorem}

It is worth noting that in \cite{RipsSela1} a more restrictive class of automorphisms is used to shorten the homomorphisms.  Since it {\em is} a more restrictive class, their results still hold using our definition of $\text{Mod}(G)$.

\begin{proposition} \label{ShortenIET}
Suppose that $Y$ is an IET subtree of $T$ and that $p_E \in \P$ is a line in $T$.  Then the intersection $Y \cap p_E$ contains at most a point.
\end{proposition}
\begin{proof}
Since $Y$ is an IET subtree, if $\sigma$ is a nondegenerate arc in $Y$ and $\epsilon > 0$ then there exists $\gamma \in \text{Stab}(Y)$ so that $\gamma . \sigma \cap \sigma$ has positive length and there is some $x \in \sigma$ such that $d_T(x, \gamma .x) < \epsilon$.  

Suppose that $Y \cap p_E$ contains more than a point.  By the above remark there exists $\gamma \in \text{Stab}(Y)$ for which $\gamma . p_E \cap p_E$ contains more than a point.  Hence $\gamma . p_E = p_E$, and $p_E \subset Y$.  This, combined with the above fact about IET components, implies that the action of $\text{Stab}(p_E)$ on $p_E$ is indiscrete.  However, this implies that it contains a noncyclic free abelian group, which cannot be a subgroup of $\text{Stab}(Y)$ when $Y$ is an IET subtree.  This contradiction proves the proposition.
\end{proof}

\begin{corollary} \label{IETinCinf}
Let $T$ be an $\R$-tree arising from some $\mathcal C_\infty$ as above.  Suppose $Y$ is an IET subtree of $T$ and $\sigma \subset Y$ a nondegenerate segment.  Then there is a segment $\hat{\sigma} \subset \mathcal C_\infty$, of the same length as $\sigma$, which corresponds to $\sigma$ under the projection from $\mathcal C_\infty$ to $T$.
\end{corollary}

\subsection{Non-IET subtrees, technical results, and the proof of Theorem \ref{ShorteningArgument}}

\ \par

An entirely analogous argument to the one which proved Proposition \ref{ShortenIET} proves

\begin{proposition}
Suppose that a line $l \subset T$ is an axial subtree and the line $p_E \subset T$ is associated to a flat $E \subset \mathcal C_\infty$.  If $l \cap p_E$ contains more than a point then $l = p_E$.
\end{proposition}

\begin{corollary} \label{NotPAxialinCinf}
Let $T$ be an $\R$-tree arising from some $\mathcal C_\infty$ as above.  Suppose $l$ is an axial component of $T$ so that $l \not\in \P$ and $\sigma \subset l$ is a non-degenerate segment.  Then there is a segment $\hat\sigma \subset \mathcal C_\infty$, of the same length as $\sigma$, which corresponds to $\sigma$ under the projection from $\mathcal C_\infty$ to $T$.
\end{corollary}

\begin{lemma} \label{DiscreteEdgeinpE}
If an edge $e$ in the discrete part of $T$ has an intersection of positive length with some line $p_E$ then $e \subset p_E$.
\end{lemma}
\begin{proof}
Suppose that $e$ contains a nontrivial segment from $p_E$ but that $e \not\subset p_E$.  Let $C$ be the edge stabiliser of $e$.  Since $\Gamma$ is freely indecomposable, $C$ is non-trivial, and since $\Gamma$ is torsion-free, $C$ is infinite.  Let $\gamma \in C$.  Then $\gamma$ leaves more than one point of $p_E$ invariant, so leaves all of $p_E$ invariant.  Thus $\gamma$ leaves $E \subset \mathcal C_\infty$ invariant.  Also, since $e \not\subset p_E$, $\gamma$ leaves some point $v \in \mathcal C_\infty \smallsetminus E$ invariant.

By \cite[Lemma 3.18]{CWIF1}, if $\{ E_i \}$ is a sequence of flats ($E_i \subset X_i$) which converges to $E$, then for all but finitely many $n$ the element $h_n(\gamma)$ leaves $E_n$ invariant.  By choosing an $n$ large enough, $h_n(\gamma).E_n = E_n$, and furthermore if $\{ v_i \}$ represents $v$, then $h_n(\gamma)$ moves $v_n$ a distance which is much smaller than the distance from $v_n$ to $E_n$.  In particular, we can ensure that the geodesic $[v_n, h_n(\gamma).v_n]$ does not intersect the $4\delta$-neighbourhood of $E_n$.  Then by \cite[Proposition 2.23]{CWIF1}, if $\pi : X_n \to E_n$ is the projection map then $d_{X}(\pi (v_n), \pi(h_n(\gamma).v_n)) \le 2\phi(3\delta)$.  However, since $h_n(\gamma) . E_n = E_n$, we know that $\pi(h_n(\gamma) . v_n) = h_n(\gamma). \pi(v_n)$. Thus $h_n(\gamma)$ moves $\pi(v_n)$ a distance at most $2\phi(3\delta)$.

Repeating this argument with a large enough subset of $C$ (namely a subset larger than the maximal size of an intersection of any orbit $\Gamma . u$ with a ball of radius $2\phi(3\delta)$), we obtain a (finite) bound on the size of $C$.  However, $C$ is infinite, as noted above.  This contradiction finishes the proof.
\end{proof}

The following Theorems \ref{AxialTheorem}, \ref{pETheorem} and \ref{DiscretewithpE} are the technical results needed to prove Theorem \ref{ShorteningArgument}.

\medskip

{\bf Theorem \ref{AxialTheorem}.}
{\em Let $G$ be a finitely generated freely indecomposable group and assume that $G \times T \to T$ is a small stable action of $G$ on an $\R$-tree $T$ with trivial stabilisers of tripods.  Let $U$ be a finite subset of $G$ and let $y \in T$.  Then there exists $\phi_A \in \text{Mod}(G)$ so that for any $u \in U$, if $[y, u.y]$ has an intersection of positive length with some axial component of $T$ then:
\[	d_T(y, \phi_A(u) . y) < d_T(y, u.y) ,	\]
and otherwise $\phi_A(u) = u$.
}

\medskip

As far as I am aware, Theorem \ref{AxialTheorem} has not appeared in print.  However, its statement and proof are very similar to those of Theorem \ref{IETTheorem}, and it is certainly known at least to Sela (see \cite[\S 5]{Sela1}) and to Bestvina and Feighn (see \cite[Exercise 11]{BFSela}).

\begin{remark} \label{fgVSfp}
Theorem \ref{IETTheorem} is stated for finitely presented groups.  The only time in the proof when it is required that $G$ be finitely presented rather than just finitely generated is when a result of Morgan from \cite{Morgan} is quoted.

Specifically, Rips and Sela show that when $G$ is freely indecomposable and finitely presented and acts on an $\R$-tree $T$ with trivial tripod stabilisers then, for each $g \in G$ and any $y \in T$, the path $[y,\gamma . y]$ cuts only finitely many components of axial or IET type in $Y$ (see \cite[pp. 350--351]{RipsSela1}).  

However, this is also true when $G$ is only assumed to be finitely generated, rather than finitely presented (but all other assumptions apply).  This follows from the arguments in \cite[\S 3]{SelaAcyl}.  The action of $G$ on $T$ can be approximated by actions of finitely presented groups $G_i$ on $\R$-trees $Y_i$.  For large enough $k$, the IET and axial components of $Y_k$ map isometrically onto the IET and axial components of $T$ (see \cite[\S 3]{SelaAcyl} for details).

Therefore, Theorem \ref{IETTheorem} stills holds when $G$ is assumed to be finitely generated, but not necessarily finitely presented.  Similarly, Theorem \ref{AxialTheorem} above, whose proof mimics the proof of Theorem \ref{IETTheorem}, holds finitely generated groups $G$.  However, in this paper we can assume $G$ is finitely presented.
\end{remark}

We now state the further technical results which are required for the proof of Theorem \ref{ShorteningArgument}.  These technical results are proved in the subsequent two sections.

\medskip

{\bf Theorem \ref{pETheorem}.}
{\em Let $\Gamma$ be a freely indecomposable, torsion-free, non-abelian toral CWIF group.  Suppose that $h_n : \Gamma \to \Gamma$ is a sequence of automorphisms converging to a faithful action of $\Gamma$ on a limiting space $\mathcal C_\infty$ and let $T$ be the $\R$-tree associated to $\mathcal C_\infty$.  Let $U$ be a finite subset of $\Gamma$, let $y \in T$ and suppose that the line $p_E \in \P$ is an axial component of $T$.  There exists a sequence of points $\{ \hat{y}_m \}$ from the $X_i$ which represents a point $\hat{y} \in \mathcal C_\infty$ and projects to $y \in T$, and there exists $m_0$ so that: for all $m \ge m_0$ there is $\phi_{p_E,m} \in \text{Mod}(\Gamma)$ so that for any $u \in U$, if $[y,u.y]$ has an intersection of positive length with a line in the $\Gamma$-orbit of $p_E$ then
\[	d_{X_m}(y_m, (h_m \circ \phi_{p_E,m})(u)) . y_m) < d_{X_m}(y_m,h_m(u). y_m) ,	\]
and otherwise $\phi_{p_E,m}(u) = u$.
}

\medskip

{\bf Theorem \ref{DiscretewithpE}.}
{\em Let $\Gamma$ be a freely indecomposable torsion-free toral CWIF group.  Suppose that $h_n : \Gamma \to \Gamma$ is a sequence of automorphisms converging to a faithful action $\Gamma$ on a limiting space $\mathcal C_\infty$ with associated $\R$-tree $T$.  Suppose further that $U$ is a finite subset of $\Gamma$ and that $y \in T$.  Then there exists a sequence $\{ \hat{y}_m \}$ representing $\hat{y} \in \mathcal C_\infty$, so that $\hat{y}$ projects to $y \in T$, and there exists $m_0$ so that: for all $m \ge m_0$ there is $\phi_{D,m} \in \text{Mod}(\Gamma)$ so that for any $u \in U$ which does not fix $y$ and with $[y,u.y]$ supported only in the discrete parts of $T$ we have
\[	d_{X_m}(\hat{y}_m, (h_m \circ \phi_{D,m})(u) . \hat{y}_m) < d_{X_m}(\hat{y}_m, u .\hat{y}_m) .	\]
}

\medskip

Armed with Theorem \ref{ShortenIET}, and assuming Theorems \ref{AxialTheorem}, \ref{pETheorem} and \ref{DiscretewithpE}, we now prove Theorem \ref{ShorteningArgument}.

\begin{proof}[Proof (Theorem \ref{ShorteningArgument}).]
We have already noted that the action of $\Gamma$ on $T$ is faithful, that $T$ is not isometric to a real line, and that the stabiliser in $\Gamma$ of any tripod in $T$ is trivial.  Also, $\SK (h_n) = \{ 1 \}$.  

We suppose (by passing to a subtree if necessary) that the tree $T$ is minimal. As noted in Remark \ref{NoThin} above, $T$ contains no thin components.

Let $U = \A_1$ be the fixed generating set of $\Gamma$ used to define $\| f \|$ for a homomorphism $f : \Gamma \to \Gamma$, and let $y$ be the image in $T$ of the basepoint $x_\omega \in \mathcal C_\infty$.

Let $\phi_I$ be the automorphism of $\Gamma$ given by Theorem \ref{IETTheorem} and $\phi_A$ the automorphism from Theorem \ref{AxialTheorem}.  

Suppose that $u \in U$ is such that $[y,u.y]$ has an intersection of positive length with an IET component of $T$.  Then Theorem \ref{IETTheorem} and Corollary \ref{IETinCinf} guarantee that for all but finitely many $n$ we have $\| h_n \circ \phi_I \| < \| h_n \|$ so $h_n$ is not short.  Similarly, if $[y,u.y]$ has an intersection of positive length with an axial component which is not contained in any $p_E \in \mathcal P$ then Theorem \ref{AxialTheorem} and Corollary \ref{NotPAxialinCinf} imply that for all but finitely many $n$ we have $\| h_n \circ \phi_A \| < \| h_n \|$ so also in this case $h_n$ is not short.

Suppose then that $[y,u.y]$ has an intersection of positive length with a line in the $\Gamma$-orbit of some $p_E$, and suppose that $p_E$ is an axial component of $T$.  Then by Theorem \ref{pETheorem} for all but finitely many $n$ there exists an automorphism $\phi_{p_E,n} \in \text{Mod}(\Gamma)$ so that $\| h_n \circ \phi_{p_E,n} \| < \| h_n \|$, so $h_n$ is not short.

Finally, suppose that all of the segments $[y, u.y]$ are entirely contained in the discrete part of $T$.  Then by Theorem \ref{DiscretewithpE} for all but finitely many $n$ there exists $\phi_{D,n} \in \text{Mod}(\Gamma)$ so that $\| h_n \circ \phi_{D,n} \| < \| h_n \|$, and once again $h_n$ is not short.

This completes the proof of the theorem.
\end{proof}

Having proved Theorem \ref{ShorteningArgument} we now prove Theorem \ref{ModinAut}.  Given Theorem \ref{ShorteningArgument}, the proof is identical to that of \cite[Corollary 4.4]{RipsSela1}.

\medskip

{\bf Theorem \ref{ModinAut}}
{\em Suppose that $\Gamma$ is a freely indecomposable torsion-free toral CWIF group.  Then $\text{Mod}(\Gamma)$ has finite index in $\text{Aut}(\Gamma)$.}

\smallskip

\begin{proof}
If $\Gamma$ is abelian then the theorem is clear, since in this case $\text{Mod}(\Z^n) = \text{Aut}(\Gamma) = \text{GL}(n,\Z)$.  Thus we assume that $\Gamma$ is non-abelian.  

For each coset $C_i = \rho_i \text{Mod}(\Gamma)$ of $\text{Mod}(\Gamma)$ in $\text{Aut}(\Gamma)$ choose a representative $\hat{\rho}_i$ which is shortest amongst all representatives of $C_i$.  That is to say, each of the automorphisms $\hat{\rho}_i$ is short. 

However, by Theorem \ref{ShorteningArgument} we cannot have an infinite sequence $\{ \hat\rho_n : \Gamma \to \Gamma \}$ of non-equivalent short automorphisms, since then some subsequence will converge to a faithful action of $\Gamma$ on a space $\mathcal C_\infty$.  Hence $\text{Mod}(\Gamma)$ has finite index in $\text{Aut}(\Gamma)$ as required.
\end{proof}

\section{Axial components} \label{AxialSection}

The purpose of this section is to prove the following two theorems.

\begin{theorem} \label{AxialTheorem}
Let $G$ be a finitely generated freely indecomposable group and assume that $G \times T \to T$ is a small stable action of $G$ on an $\R$-tree $T$ with trivial stabilisers of tripods.  Let $U$ be a finite subset of $G$ and let $y \in T$.  Then there exists $\phi_A \in \text{Mod}(G)$ so that for any $u \in U$, if $[y, u.y]$ has an intersection of positive length with some axial component of $T$ then:
\[	d_T(y, \phi_A(u) . y) < d_T(y, u.y) ,	\]
and otherwise $\phi_A(u) = u$.
\end{theorem}

\begin{theorem} \label{pETheorem}
Let $\Gamma$ be a freely indecomposable, torsion-free, non-abelian toral CWIF group.  Suppose that $h_n : \Gamma \to \Gamma$ is a sequence of automorphisms converging to a faithful action of $\Gamma$ on a limiting space $\mathcal C_\infty$ and let $T$ be the $\R$-tree associated to $\mathcal C_\infty$.  Let $U$ be a finite subset of $\Gamma$, let $y \in T$ and suppose that the line $p_E \in \P$ is an axial component of $T$.  There exists a sequence of points $\{ \hat{y}_m \}$ from the $X_i$ which represents a point $\hat{y} \in \mathcal C_\infty$ and projects to $y \in T$, and there exists $m_0$ so that: for all $m \ge m_0$ there is $\phi_{p_E,m} \in \text{Mod}(\Gamma)$ so that for any $u \in U$, if $[y,u.y]$ has an intersection of positive length with a line in the $\Gamma$-orbit of $p_E$ then
\[	d_{X_m}(y_m, (h_m \circ \phi_{p_E,m})(u)) . y_m) < d_{X_m}(y_m,h_m(u). y_m) ,	\]
and otherwise $\phi_{p_E,m}(u) = u$.
\end{theorem}

To prove Theorem \ref{AxialTheorem} we follow the proof of \cite[Theorem 5.1]{RipsSela1} (which is Theorem \ref{IETTheorem} in this paper).  First, we need the following result, the (elementary) proof of which we include because of its similarity to Proposition \ref{ShortenpEAxial} below.

\begin{proposition} \label{ShortenAxial}
Suppose that $\rho: P \times \R \to \R$ is an orientation-preserving, indiscrete isometric action of $P \cong \Z^n$ on the real line $\R$.  For any finite subset $W$ of $P$and any $\epsilon > 0$ there exists an automorphism $\sigma : P \to P$ such that:
\begin{enumerate}
\item[1)] For every $w \in W$ and every $r \in \R$
\[	d_{\R}(r, \sigma (w) . r ) < \epsilon ;	\]
\item[2)] For any $k \in \text{ker}(\rho)$ we have $\sigma(k) = k$.
\end{enumerate}
\end{proposition}
\begin{proof}
There is a direct product decomposition $P = A \oplus B$ where $A$ is the kernel of the action of $P$ on $\R$, and $B$ is a finitely generated free abelian group which has a free, indiscrete and orientation preserving action on $\R$.  The automorphism $\sigma$ we define fixes $A$ elementwise, so we can assume that all elements of $W$ lie in $B$ (since elements of $A$ fix $\R$ pointwise).  Thus, we need only prove the lemma in case the action is faithful.

Since the action of $B$ on $\R$ is indiscrete and free, the translation lengths of elements of a basis of $B$ are independent over $\Z$.  In particular, there is a longest translation length amongst the translation lengths of a basis of $B$.  Suppose that $b_1 \in B$ is the element of the basis with largest translation length, and that $b_2$ has the second largest.  Denote these translation lengths by $| b_1 |$ and $| b_2 |$, respectively.  Since $| b_1 |$ and $| b_2 |$ are independent over $\Z$, there is $n \in \Z$ so that $0 < |b_1 + nb_2| < | b_2 |$.  Replace $b_1$ by $b_1 + nb_2$.  This is an automorphism of $P$, fixing $A$ elementwise.  

Proceeding in this manner, we can make each of the elements of a basis as small as we wish, and so given $W$ and $\epsilon > 0$, we can make each of the elements of $W$ (considered as a word in the basis of $B$) have translation length less than $\epsilon$, as required.
\end{proof}

\begin{proof}[Proof (Theorem \ref{AxialTheorem}).]

By Claim \ref{fgVSfp}, each of the segments $[y,u.y]$ for $u \in U$ cuts only finitely many components of $T$ of axial or IET type.  Let $\epsilon$ be the minimum length of a (non-degenerate) interval of intersection between $[y,u.y]$ and an axial component of $T$, for all $u \in U$.

The action of $G$ on $T$ induces a graph of groups decomposition $\Lambda$ of $G$ as in Theorem \ref{RtreeStructure}.  Let $T_i$ be an axial component of $T$.  There is a vertex group of $\Lambda$ corresponding to the $G$-orbit of $T_i$, with vertex group a conjugate of $\text{Stab}(T_i)$.  By Theorem \ref{Texists} and Lemma \ref{LimitCSA}, $\text{Stab}(T_i)$ is a free abelian group.  The vertex groups adjacent to $\text{Stab}(T_i)$ (those that are separated by a single edge) stabilise a point in the orbit of a branching point in $T_i$ with nontrivial stabiliser.  Recall that $G$ is freely indecomposable, so all edge groups are nontrivial.

Let $q_1$ be the point on $T_i$ closest to $y$ (if $y \in T_i$ then $q_1 = y$).  Choose points $q_2, \ldots , q_m \in T_i$ in the orbits of the branching points corresponding to the adjacent vertex groups such that $d_T(q_i,q_j) < \frac{\epsilon}{20}$.  We can do this since the action of $\text{Stab}(T_i)$ on $T_i$ has all orbits dense, since $T_i$ is an axial component.

The proof of Theorem \ref{AxialTheorem} now proceeds exactly as the proof of \cite[Theorem 5.1]{RipsSela1} (start with Case 1 on p.351).
\end{proof}

\begin{proof} [Proof (Theorem \ref{pETheorem}).]
Since $p_E \in \P$ is an axial component of $T$, there is a vertex group corresponding to the conjugacy class of $\text{Stab}(p_E)$ in the graph of groups decomposition which the (faithful) action of $\Gamma$ on $T$ induces (see Theorem \ref{RtreeStructure}).

Now, the stabiliser in $\Gamma$ of $p_E$ is exactly the stabiliser in $\Gamma$ of $E$, when $\Gamma$ acts (also faithfully) on $\mathcal C_\infty$.  By \cite[Corollary 3.17]{CWIF1}, there is a sequence of flats $E_i$ in the approximating spaces $X_i$ so that $E_i \to E$ in the $\Gamma$-equivariant Gromov topology. By \cite[Lemma 3.18]{CWIF1}, if $\gamma \in \text{Stab}_\Gamma(E)$ then for all but finitely $i$ we have $h_i(\gamma) \in \text{Stab}(E_i)$. For such an $i$, the element $h_i(\gamma)$ is contained in a unique noncyclic maximal abelian subgroup $A_i$  of $\Gamma$.  
However, $h_i$ is an automorphism, so $\gamma$ is contained in a unique noncyclic maximal abelian subgroup $A_\gamma$ of $\Gamma$, and $A_i = h_i(A_\gamma)$.

If $\gamma'$ is another element of $\text{Stab}_\Gamma(E)$, then $[\gamma, \gamma'] = 1$, and it is not difficult to see that $A_\gamma = A_{\gamma'}$.  Also, if $\gamma_0 \in A_\gamma$ then $\gamma_0 \in \text{Stab}_\Gamma(E)$.  Hence $A_\gamma = \text{Stab}_\Gamma(E)$. We denote the subgroup $\text{Stab}_\Gamma(E)$ by $A_E$.

We now prove Theorem \ref{pETheorem} by finding an analogue of Proposition \ref{ShortenAxial} in the flats $E_i$ and then once again following the proof from \cite{RipsSela1}.

\begin{proposition} \label{ShortenpEAxial}
Let $W$ be a finite subset of $A_E$.  For any $\epsilon > 0$ there exists $i_0$ so that for all $i \ge i_0$, there is an automorphism $\sigma_i : A_E \to A_E$ so that 
\begin{enumerate}
\item For every $w \in W$, and every $r_i \in E_i$,
\[	d_{X_i}(r_i, h_i(\sigma(w)) . r_i) < \epsilon ;	\]
\item For any $k \in A_E$ which acts trivially on $E$ we have $\sigma(k) = k$.
\end{enumerate}
\end{proposition}
\begin{proof}[Proof (Proposition \ref{ShortenpEAxial}).]
The group $A_E$ admits a decomposition $A_E = A_0 \oplus A_1$, where $A_0$ acts trivially on $E$, and $A_1$ acts freely on $E$.  Choose a basis $\B$ of $A_E$ consisting of a basis for $A_0$ and a basis for $A_1$.  Let $k_W$ be the maximum word length of any element of $W$ with respect to the chosen basis.

Since the $h_i : \Gamma \to \Gamma$ are automorphisms, for sufficiently large $i$ and any $a \in E_i$, the set $h_i(A_E) . a \subset E_i$ forms an $\frac{\epsilon}{20k_W}$-net in $E_i$ (where distance is measured in the metric $\frac{1}{\| h_i \|}$ on $X_i$).  Choose a (possibly larger) $i$ so that also the action of $h_i(\B)$ on $E_i$ approximates the action of $\B$ on $E$ to within $\frac{\epsilon}{20k_W}$ (note that since the action of $A_E$ on $E$ and the action of $h_i(A_E)$ on $E_i$ are both by translations, and translations of Euclidean space move every point the same distance, there are arbitrarily good approximations for the action of any finite subset of $A_E$ on the whole of $E$).

The remainder of the proof proceeds just as the proof of Propostion \ref{ShortenAxial} above, although in the step where we replace $b_1$ by $b_1 + nb_2$, we cannot insist that $b_2$ acts nontrivially on $E$.  However, we of course can insist that $b_1$ acts nontrivially on $E$, since otherwise it moves all points of $E$ a distance at most $\frac{\epsilon}{20}$.  Therefore, such an automorphism is nonetheless a generalised Dehn twist.
\end{proof}
Given Proposition \ref{ShortenpEAxial}, the proof of Theorem \ref{pETheorem} once again follows the proof of \cite[Theorem 5.1, pp. 350--353]{RipsSela1}, although in this case we have to choose approximations to the action of $\Gamma$ on $\mathcal C_\infty$ (the important point here is that the sets $h_i(A_E) . a$, for any $a \in E_i$, get denser and denser in $E_i$, when considered in the scaled metric $\frac{1}{\| h_i \|}d_X$ of $X_i$).  These small changes are straightforward, but do lead to the different shortening automorphisms $\phi_{p_E,m}$ in the statement of Theorem \ref{pETheorem}.
\end{proof}

\section{The discrete case} \label{Discrete}

In this section we shorten the approximations to paths of the form $[\hat{y},u. \hat{y}]$, where $\hat{y} \in \mathcal C_\infty$ projects to $y \in T$ and $[y,u.y]$ is entirely supported in the discrete part of $T$.  The lengths of the limiting paths $[\hat{y}, u. \hat{y}]$ and $[y,u.y]$ are unchanged.

The purpose of this section is to prove the following

\begin{theorem} \label{DiscretewithpE}
Let $\Gamma$ be a freely indecomposable torsion-free toral CWIF group.  Suppose that $h_n : \Gamma \to \Gamma$ is a sequence of automorphisms converging to a faithful action $\Gamma$ on a limiting space $\mathcal C_\infty$ with associated $\R$-tree $T$.  Suppose further that $U$ is a finite subset of $\Gamma$ and that $y \in T$.  Then there exists a sequence $\{ \hat{y}_m \}$ representing $\hat{y} \in \mathcal C_\infty$, so that $\hat{y}$ projects to $y \in T$, and there exists $m_0$ so that: for all $m \ge m_0$ there is $\phi_{D,m} \in \text{Mod}(\Gamma)$ so that for any $u \in U$ which does not fix $y$ and with $[y,u.y]$ supported only in the discrete parts of $T$ we have
\[	d_{X_m}(\hat{y}_m, (h_m \circ \phi_{D,m})(u) . \hat{y}_m) < d_{X_m}(\hat{y}_m, u .\hat{y}_m) .	\]
\end{theorem}

The proof of Theorem \ref{DiscretewithpE} follows \cite[\S 6]{RipsSela1}.

By Lemma \ref{DiscreteEdgeinpE}, if $e$ is a discrete edge in $T$ then either $e \in p_E$ for some flat $E \subset \mathcal C_\infty$, or there is a well-defined, canonical, isometric image $\hat{e}$ of $e$ in $\mathcal C_\infty$, so that $\hat{e}$ projects to $e$.

We have a sequence of automorphisms $\{ h_n : \Gamma \to \Gamma \}$, converging to a faithful action of $\Gamma$ on a limiting space $\mathcal C_\infty$, with associated $\R$-tree $T$. 

Since the length of a homomorphism is measured only by its action on the discrete set $\Gamma . x \subset X$, it is worth noting that in all of the analysis below, it is possible to choose a sequence $\{ \hat{y}_m \}$ which represents $\hat{y}$, so that $\hat{y}_m \in \Gamma . x$ for all $m$ (and $\hat{y} \in \mathcal C_\infty$ projects to $y \in T$).  Thus we can ensure that our `shortening' automorphisms really do shorten.  We mostly do this without comment (and we implicitly used this in the proof of Theorem \ref{ShorteningArgument}).

There are a number of different cases to consider:

\smallskip

\noindent{\bf \underline{Case 1:}}  $y$ is contained in the interior of an edge $e$

\smallskip

\noindent{\bf Case 1a:}  $e$ is not completely contained in a line of the form $p_E$ and $\bar{e} \in T/\Gamma$ is a splitting edge.  

Note that because $e$ is not contained in any $p_E$, there is a single point $\hat{y} \in \mathcal C_\infty$ which corresponds to $y \in T$.

This case is very similar to the Case 1a on pp. 355--356 of \cite{RipsSela1}.  In this case we have a decomposition $\Gamma = A \ast_C B$ where $C$ is a finitely generated free abelian group properly contained in both $A$ and $B$.

Given $u \in U$ we can write:
\[	u = a_u^1b_u^1 \cdots a_u^{n_u}b_u^{n_u} ,	\]
where $a_u^i \in A$ and $b_u^i \in B$ (possibly $a_u^1$ and/or $b_u^{n_u}$ are the identity).  Let $\{ z_1, \ldots , z_n \}$ be a generating set for $Z$.

Let $\epsilon$ be the minimum of:
\begin{enumerate}
\item the length of the shortest edge in the discrete part of $T$;
\item the distance between $y$ and the vertices of $e$.
\end{enumerate}

Recall that triangles in $X$ are relatively $\delta$-thin, and the function $\phi$ comes from the definition of isolated flats.  Let $C_0$ be the maximum size of an intersection of an orbit $\Gamma . z$ with a ball of radius $10 \delta + 2\phi(3\delta)$ in $X$ (where distance is measured in $d_X$).

Now take $F$ to be the finite subset of $G$ containing $1$ and
\[	z a_u^i z^{-1},\  \ \ zb_u^iz^{-1} ,	\]
where $z \in C$ has word length at most $10C_0$.

For large enough $m$ we have, for all $\gamma_1 , \gamma_2 \in F$,

\begin{equation} \label{Approx}
| d_{X_m}(h_m(\gamma_1) . y_m, h_m(\gamma_2).y_m) - d_T(\gamma_1 . y, \gamma_2 . y)| < \epsilon_1 ,
\end{equation}
where $\epsilon_1 = \frac{\epsilon}{8000C_0}$.

Let $w_m \in [y_m,h_m(a_u^1).y_m]$ and $w_m' \in [y_m, h_m(b_u^1).y_m]$ satisfy 
\begin{equation} \label{Flexible}
\frac{\epsilon}{2} - \frac{\epsilon}{1000} \le d_{X_m}(w_m,y_m) = d_{X_m}(w_m',y_m) \le \frac{\epsilon}{2} + \frac{\epsilon}{1000}.
\end{equation}

\begin{lemma} \label{zMoves}
For some $z \in C$ of word length at most $10C_0$ we have, for all but finitely many $m$,
\begin{eqnarray*}
d_{X_m}(y_m, h_m(z). w_m) & < & d_{X_m}(y_m,w_m) - 8\delta_m, \mbox{ and } \\
d_{X_m}(y_m, h_m(z).w_m') & < & d_{X_m}(y_m,w_m')+ 8\delta_m.
\end{eqnarray*}
\end{lemma}
\begin{proof}
Let $W$ be the set of all elements $z \in C$ of word length at most $10C_0$ in the generators $\{ z_1, \ldots , z_n \}$ and their inverses.

First suppose that for all but finitely many $i$ we have $h_i(W) \subseteq \text{Stab}_\Gamma(E_i)$.  Then since the edge containing $y$ is not completely contained in a single $p_E$, we can assume that each element of $W$ fixes a point outside of $E$.  Now, using \cite[Proposition 2.23]{CWIF1}, there is a point in $E_i$ which is moved at most $2\phi(3\delta)$ by each element of $h_i(W)$.  This gives a bound on the size of $h_i(W)$ which does not depend on $i$ (so long as $i$ is large enough).  However, this contradicts the choice of $W \subseteq \Gamma$.  Therefore it is not the case that $h_i(W) \subseteq \text{Stab}_\Gamma(E_i)$ for all but finitely many $i$.

By the argument in the paragraphs after the proof of Lemma 4.5, for all but finitely many $k$, the elements $h_i(z)$ act approximately like translations.  Since $W$ is closed under inverses, and we have chosen $W$ large enough that some element `translates' by at least $10\delta_m$, we can choose some $z \in W$ which satisfies the conclusion of the lemma.
 \end{proof}
 
 In order to finish Case 1a, we follow the proofs of Proposition 6.3 and Theorem 6.4 from \cite{RipsSela1}.  The only additional thing needed in this case is to force $w_m$ to lie close to each $[y_m, h_m(a_u^i).y_m]$.  We do this by applying \cite[Lemma 4.5]{CWIF1} and the arguments in the paragraphs in \cite{CWIF1} which follow the proof of that result.  It is for this reason that we left some flexibility as to the choice of $w_m$ and $w_m'$ in \ref{Flexible} above.  
 
 Doing this, we can ensure that $w_m$ lies within $2\delta$ of each geodesic segment $[y_m, h_m(a_u^i) . y_m]$, and similarly for $w_m'$.  We can now follow the proof \cite[Proposition 6.3]{RipsSela1}.  The proof of \cite[Theorem 6.4]{RipsSela1} is not included in \cite{RipsSela1} (or in \cite{SelaAcyl} as claimed in \cite{RipsSela1}).  However, it is straightforward, so we omit it here also.
 
The automorphism we use to shorten in this case is:
  \begin{eqnarray*}
 \forall a \in A && \phi(a) = zaz^{-1} \\
 \forall b \in B && \phi(b) = z^{-1}bz,
 \end{eqnarray*}
 where $a$ is as in Lemma \ref{zMoves} above.  This completes the proof in Case 1a.  It is worth noting here that we are shortening the actions on $X_i$ which approximate the action on $\mathcal C_\infty$.  However, this does not affect the analogy between the proofs here and those in \cite{RipsSela1}.

\smallskip

\noindent{\bf Case 1b:} $e$ is not completely contained in a single $p_E$ and $\bar{e} \in T/\Gamma$ is not a splitting edge.  

In this case we have a decomposition $\Gamma = A \ast_C$, where $C$ is a finitely generated free abelian group.  

In the same way as we adapted the proof of Case 1a from \cite{RipsSela1} above, we may adapt the proof of Case 1b from \cite{RipsSela1}.  The key point is that we allow a small amount of flexibility in the choice of $w_m$ and $w_m'$.  Doing this, we may ensure that even though the approximating triangles we consider are only relatively thin, rather than actually thin, all of the features we need to apply the proof from \cite{RipsSela1} still hold, because we can make sure that we are not near the `fat' part of any triangle.  Proceeding with this idea in mind, the proof from \cite{RipsSela1} can be adapted without difficulty.

\medskip

\medskip

We now deal with the two cases where $y$ is contained in the interior of the edge $e$ and $e \subset p_E$ for some $p_E \in \P$.  Using \cite[Lemma 2.22]{CWIF1} and \cite[Proposition 2.23]{CWIF1}, the following result is not difficult to prove:

\begin{proposition}
Suppose that $X$ is a CWIF space, and that $E_1, E_2 \in \mathcal F$ are maximal flats in $X$.  There exists a bounded set $J_{E_1,E_2}$ so that for any $x \in E_1$ and any $y \in E_2$, the geodesic $[x,y]$ intersects $J_{E_1,E_2}$.  Moreover, we can choose $J_{E_1,E_2}$ so that
\[	\text{Diam}(J_{E_1,E_2}) \le 7\delta + 14\phi(4\delta)	.	\]
\end{proposition}
Recall that $\mathcal F$ is the family of maximal flats from the definition of isolated flats, that triangles in $X$ are relatively $\delta$-thin, and that $\phi$ is the function from the definition of isolated flats.  We assume without loss of generality that for all $k \ge 0$ we have $\phi(k) \ge k$ and also that $\phi$ is a nondecreasing function.

Choose compact fundamental domains for the action of $\text{Stab}_\Gamma(E)$ on $E$, for each conjugacy class of maximal flat in $X$, and let $K_F$ be the maximal diameter of these fundamental domains.  Also, let $K_X$ be the diameter of a compact set $D$ for which $\Gamma . D = X$.  For the remainder of Case 1, we replace the constant $\delta$ by
\[	\max \left\{ \delta , 1000K_F , 1000 K_X, 1000(7\delta + 14\phi(4\delta)) \right\}	.	\]

The stabiliser of the edge $e$ is a subgroup of $\text{Stab}_\Gamma(E)$.  Since $p_E$ is not an axial component, the action of $\text{Stab}(E)$ on $E$ is either trivial or factors through a infinite cyclic group.  If $\bar{e} \in T/\Gamma$ is a splitting edge, then necessarily the action of $\text{Stab}(E)$ on $E$ is trivial. 

\noindent{\bf Case 1c:}  $e$ is completely contained in some $p_E$, and $\bar{e} \in T/\Gamma$  is a splitting edge.

Let $A_E = \text{Stab}_\Gamma(E)$.  Then, $A_E = A_0 \oplus A_1$, where $A_0$ acts trivially on $E$ and $A_1$ acts freely on $E$.  Since $p_E$ is a splitting edge, $A_1 = \{ 1 \}$.

We have a decomposition $\Gamma = H_1 \ast_{A_E} H_2$.

The subgroup $H_1$ fixes a point in $p_E$, but does not fix all of $p_E$.  Thus, $H_1$ fixes a point $v_1 \in E$.  Similarly, $H_2$ fixes a point $v_2 \in E$, but does not fix all of $E$.  We choose points $\hat{y}_m \subset E_m$ so that: (i) $\{ \hat{y}_m \}$ represents $\hat{y} \in \mathcal C_\infty$ which projects to $y \in T$; (ii) each $\hat{y}_m$ lies in the orbit $\Gamma . x$; and (iii) subject to the first two conditions, $\hat{y}_m$ lies as close as possible to the line $[v_1^m, v_2^m]$, where $\{ v_i^m \} \to v_i$, $i = 1,2$.

We proceed as in Case 1a. However, this time we cannot find a single automorphism to shorten the $\| h_i \|$, but we use the fact that the sets $h_i(A_E) .a \subset E_i$ are denser and denser (when distance is measured in the metrics $\frac{1}{\| h_i \|}d_X$) to find, for all but finitely many $i$, a Dehn twist $\phi_{e,i}$ which shortens the action on $X_i$.  This proceeds in a similar way to Case 1a above, using the ideas in Proposition \ref{ShortenpEAxial} and the proof of Theorem \ref{pETheorem} above.

\smallskip

\noindent{\bf Case 1d:} $e$ is completely contained in some $p_E$ and $\bar{e} \in T/\Gamma$ is not a splitting edge.

There are two cases here.  As in Case 1b, we have a decomposition $\Gamma = A \ast_C$, where $C$ is a finitely generated free abelian group.  Let $t$ be the stable letter of this HNN extension, and suppose that $C \le \text{Stab}(E)$, a maximal flat in $\mathcal C_\infty$.  The two cases are where $f \in \text{Stab}(E)$, and when $f \not\in \text{Stab}(E)$.

Each of these cases follow the proof of Case 1b above (and therefore Case 1b from \cite{RipsSela1}) in the same way as Case 1c followed the proof of Case 1a.

\medskip

\noindent{\bf Case 2:}  $y$ is a vertex of $T$.

In this case, we do not shorten the approximations to a particular edge, but each of the edges adjacent to $y$.  As before, we largely follow \cite[\S 6]{RipsSela1}.

There are four cases again, when the edge is splitting, and non-splitting, coupled with the cases where the edge is contained in some $p_E$ and when it is not.

These follow the proofs from \cite{RipsSela1} just as in Case 1 above.  Note that the shortening automorphisms fix elementswise any element which fixes the point $\hat{y} \in \mathcal C_\infty$ which projects to $y \in T$.

\begin{proof}[Proof (Theorem \ref{DiscretewithpE})]
If $y$ is contained in the interior of an edge, then apply Case 1 above to find a sequence of automorphisms which shorten the $h_n$.

If $y$ is a vertex in $T$, then we shorten the $h_n$ on each of the adjacent edges separately using Case 2 and \cite[\S 6]{RipsSela1}.
\end{proof}

This finally finishes the proof of Theorem \ref{ShorteningArgument}.

\section{The Hopf property} \label{HopfSection}

In this section we prove the main result of this paper:

\medskip

{\bf Theorem \ref{GammaHopfian}.}
{\em Suppose that $\Gamma$ is a torsion-free toral CWIF group.  Then $\Gamma$ is Hopfian.
}

\medskip

We prove Theorem \ref{GammaHopfian} by assuming that there is a surjective homomorphism which is not an automorphism and deriving a (rather involved) contradiction.

Suppose that $\Gamma$ is a torsion-free toral CWIF group and that $\psi : \Gamma \to \Gamma$ is a surjective homomorphism and that $\text{ker}(\psi) \neq \{ 1 \}$.  Note that finitely generated free abelian groups are Hopfian, so $\Gamma$ is nonabelian.

Let $H$ be a finitely generated subgroup of $\Gamma$, and assume that $\psi$ restricts to an epimorphism $\psi_H$ from $H$ onto itself, and that $\psi_H$ has a nontrivial kernel.  Note that for the moment $H$ could be $\Gamma$.  Other subgroups satisfying these assumptions arise later in the proof.  Once again, since $H$ is not Hopfian, $H$ is nonabelian.

The homomorphisms $\psi_H^m : H \to \Gamma$ have different kernels, so $\psi_H^n$ and $\psi_H^m$ are non-conjugate when $m \neq n$.  Thus, we can apply the construction from \cite{CWIF1} to obtain a limiting space $\mathcal C_{\infty,H}$, and an $\R$-tree $T_H$, both equipped with isometric actions of $H$ with no global fixed point.  Let $K_H$ be the kernel of the action of $H$ on $\mathcal C_{\infty,H}$ (which is also the kernel of the action of $H$ on $T_H$).

In the course of this construction, we conjugate each $\psi_H^m$ by some $\gamma_m \in \Gamma$ so that $\phi_m := \tau_{\gamma_m} \circ \psi_H^m$ has `minimal displacement'.  Since this conjugation does not change the kernel we have $\text{ker}(\phi_m) = \text{ker}(\psi_H^m) = \text{ker}(\psi^m) \cap H$.

Define $H_\infty = H / K_H$, a strict $\Gamma$-limit group, and let $\eta : H \to H_\infty$ be the natural surjection.  The following result follows immediately from Theorem \ref{Texists}.

\begin{lemma} \label{HinfProps} 
\ \par

\begin{enumerate}
\item If $h \in H$ stabilises a tripod in $T_H$ then $h \in \text{ker}(\psi^m)$ for some $m$;
\item The stabiliser in $H_\infty$ of any non-degenerate segment in $T_H$ is abelian; \label{SegStabAb}
\item Let $[y_1, y_2] \subset [y_3,y_4]$ be a pair of non-degenerate segments of $T_H$ and assume that the pointwise stabiliser $\text{Fix}_{H_\infty}([y_3,y_4])$ is non-trivial.  Then
\[	\text{Fix}_{H_\infty}([y_1,y_2]) = \text{Fix}_{H_\infty}([y_3,y_4]) .	\]
In particular, the action of $H_\infty$ on $T_H$ is stable. Finally,
\item $H_\infty$ is torsion-free.
\end{enumerate}
\end{lemma}

We collect some elementary properties about $H_\infty$ and $T_H$.

\begin{lemma} \cite[Corollary 5.7.4]{CWIF1}
The group $H_\infty$ is CSA.  In other words, if $A$ is a maximal abelian subgroup of $H_\infty$ then $A$ is malnormal.
\end{lemma}

\begin{lemma} [cf. Lemma 1.4, p.305, \cite{SelaHopf}]
The following properties hold:
\begin{enumerate}
\item $K_H = \bigcup\limits_{k=1}^\infty \text{ker}(\phi_k)$;
\item $H_\infty$ is finitely generated and $\text{rk}(H_\infty) = \text{rk}(H)$; and
\item $H_\infty$ is not finitely presented.
\end{enumerate}
\end{lemma}
\begin{proof}
Identical to the proof of \cite[Lemma 1.4, p. 305]{SelaHopf}.
\end{proof}

\begin{proposition} [cf. Proposition 1.9, p.307, \cite{SelaHopf}] \label{NoIET}
If $H_\infty$ is freely indecomposable then $T_H$ satisfies the following properties:
\begin{enumerate}
\item $T_H$ is not isometric to a real line;
\item $T_H$ does not contain any thin components; and
\item Let $\sigma$ be a non-degenerate segment in $T_H$ with a nontrivial stabiliser in $H_\infty$.  Then $\sigma$ does not have an intersection of positive length with any IET component of $T_H$.
\end{enumerate}
\end{proposition}

\begin{lemma} [cf. Lemma 1.10, p.308, \cite{SelaHopf}]
For each $x \in H_\infty$ let $h_x \in H$ satisfy $\eta (h_x) = x$.  let $\nu_\infty : H_\infty \to H_\infty$ be a map given by $\nu_\infty(x) = \eta(\psi(h_x))$.  Then:
\begin{enumerate}
\item $\nu_\infty$  is well-defined:  if $x = \eta(h_1) = \eta(h_2)$ then $\eta(\psi(h_1)) = \eta(\psi(h_2)) = \nu_\infty(x)$;
\item $\nu_\infty$ is an endomorphism of $H_\infty$; and
\item $\nu_\infty$ is an automorphism of $H_\infty$.
\end{enumerate}
\end{lemma}
\begin{proof}
Identical to the proof of \cite[Lemma 1.10, p.308]{SelaHopf}.
\end{proof}

Let $q$ be the maximal rank of a free abelian subgroup of $\Gamma$.  In addition to Lemma \ref{HinfProps}.(\ref{SegStabAb}), we have

\begin{lemma} [cf. Lemma 1.3, p.304, \cite{SelaHopf}]
Stabilisers of non-degenerate segments of $T$ in $H_\infty$ are either trivial or locally (free abelian of rank at most $q$).
\end{lemma}
\begin{proof}
Recall that  $T$ is not isometric to a real line, so $K_H = \bigcup\limits_{k=1}^\infty \text{ker}(\phi_k)$.

Let $A$ be the stabiliser in $H_\infty$ of a non-degenerate segment of $T$, let $\hat{A} \in H$ satisfy $\eta (\hat{A}) = A$ and let $a_1, \ldots, a_n \in \hat{A}$.  Then \cite[Theorem 4.4.(1)]{CWIF1} shows, for some large $k$, that $\langle \psi^k(a_i), \ldots , \psi^k(a_n) \rangle$ is an abelian subgroup of $\Gamma$, and thus is free abelian of rank at most $q$.  The result now follows from the fact that $\eta \circ \psi^k = \nu_\infty^k \circ \eta$, since $\nu_\infty$ is an automorphism.
\end{proof}

\subsection{The JSJ decomposition of $H_\infty$}

Since $H_\infty$ is a strict $\Gamma$-limit group, Theorem \ref{AbelianJSJ} holds for $H_\infty$. We would like to consider how the abelian JSJ decomposition $\Lambda_{H_\infty}$ of $H_\infty$ is affected by the automorphism $\nu_\infty$.  Certainly $\nu_\infty$ induces another graph of groups decomposition of $H_\infty$.  

We first consider the vertex groups.  Let $V_1, \ldots , V_m$ be the vertex groups of $\Lambda_{H_\infty}$.  If $V_i$ is a CMQ subgroup of $\Lambda_{H_\infty}$ then $\nu_\infty(V_i)$ is also a CMQ subgroup.  By Theorem \ref{AbelianJSJ}.(1), in this case $\nu_\infty(V_i)$ is a conjugate of some vertex group of $\Lambda_{H_\infty}$.  If $V_i$ is not a CMQ subgroup then by Theorem \ref{AbelianJSJ}.(1) $\nu_\infty(V_i)$ is elliptic in $\Lambda_{H_\infty}$, so is again conjugate to some vertex group of $\Lambda_{H_\infty}$.

Now consider the edge groups of $\Lambda_{H_\infty}$.  Suppose that $E_j$ is a non-cyclic edge group in $\Lambda_{H_\infty}$, with $V_{j_1}$ and $V_{j_2}$ the adjacent vertex groups.  Then $\nu_\infty(E_j)$ is elliptic in $\Lambda_{H_\infty}$ and is contained in $\nu_\infty(V_{j_1})$ and $\nu_\infty(V_{j_2})$, and thus $\nu_\infty(E_j)$ is a conjugate of an edge group of $\Lambda_{H_\infty}$.  If $E_j$ is a cyclic edge group in $\Lambda_{H_\infty}$ then the above argument applies unless $\nu_\infty(E_j)$ is hyperbolic in $\Lambda_{H_\infty}$.  However, in this case the generator of $\nu_{\infty}(E_j)$ corresponds to a weakly essential s.c.c. in the $2$-orbifold corresponding to some CMQ subgroup of $H_\infty$.  This cannot happen by the maximality of CMQ subgroups.  Thus we have the following:

\begin{proposition} [cf. Corollary 2.9, p.311, \cite{SelaHopf}] \label{VertexPeriodic}
Suppose that $H_\infty$ is freely indecomposable.  Then the automorphism $\nu_\infty$ permutes the conjugacy classes of vertex and edge groups in the abelian JSJ decomposition of $H_\infty$.  Furthermore, if $E$ is a maximal cyclic edge group in the abelian JSJ decomposition of $H_\infty$ and $e \in E$, then the conjugacy class of $e$ is periodic under $\nu_\infty$.
\end{proposition}

\subsection{Proving Theorem \ref{GammaHopfian}}

Let $\Gamma$ be a torsion-free toral CWIF group and suppose that $\Gamma$ is not Hopfian.  Let $\psi : \Gamma \to \Gamma$ be an epimorphism with nontrivial kernel, and let $\Gamma_\infty$ be the corresponding strict $\Gamma$-limit group as constructed above.  

Note that $\Gamma_\infty$ is finitely generated but not finitely presented, so $\eta : \Gamma \to \Gamma_\infty$ is not an isomorphism.  Also, $\Gamma_\infty$ is not a free group, so we have:
\[	\Gamma_\infty = H_\infty^1 \ast \cdots \ast H_\infty^l \ast F_r ,	\]
where $H_\infty^1, \ldots , H_\infty^l$ are finitely generated, non-cyclic, torsion-free and freely indecomposable and $F_r$ is an additional free factor (possibly trivial).  Note that $l \ge 1$ since $\Gamma_\infty$ is not a free group.

There is $m \ge 1$ and elements $\gamma_1, \ldots , \gamma_l \in \Gamma$ so that for each $1 \le i \le l$ we have $\nu_\infty^m (H_\infty^i) = \eta(\gamma_i) H_\infty^i \eta(\gamma_i^{-1})$.

Replacing $\psi$ by $\psi^m$ composed with an inner automorphism of $\Gamma$ (and keeping the same notation), we assume that $\nu_\infty$ preserves $H_\infty^1$.  The subgroup $H_\infty^1$ is finitely generated, so let $h_1, \ldots , h_t \in \Gamma$ be such that $H_\infty^1 = \langle \eta(h_1) , \ldots , \eta(h_t) \rangle$.  By the argument in the last paragraph of page 311 of \cite{SelaHopf}, there is some $k$ so that if we define $H_1 = \langle \psi^k(h_1), \ldots \psi^k(h_t) \rangle$ then $\psi(H_1) = H_1$ and $\eta(H_1) = H_\infty^1$. 

If $H_\infty^1$ is finitely presented then it is not difficult to see that $\eta |_{H_1} : H_1 \to H_\infty^1$ is an isomorphism (since $K_H = \bigcup\limits_{k=1}^\infty \text{ker}(\psi^m)$).  Therefore, we suppose for now that $H_\infty^1$ is not finitely presented.  In this case we have that (in our previous notation) $(H_1)_\infty = H_\infty^1$.

\begin{theorem} [cf. Theorem 3.1, p.312, \cite{SelaHopf}]
Let $H_1$ and $H_\infty^1$ be as above, and let $\Lambda_{H_\infty^1}$ be the abelian JSJ decomposition of $H_\infty^1$.  If all the edge groups in $H_\infty^1$ are finitely generated then $\eta |_{H_1} : H_1 \to H_\infty^1$ is an isomorphism.
\end{theorem}
\begin{proof}
Similar to the proof of \cite[Theorem 3.1, p.312]{SelaHopf}.  Note that we may need to invoke the arguments in Paragraph \ref{AcylSection} above in order to ensure that edge groups adjacent to the vertex group $V_\infty^i$ are elliptic in the graph of groups decomposition found during the proof in \cite{SelaHopf}.  This is because the finitely generated edge groups may be non-cyclic, so that the conjugacy classes of the generators need not be periodic under the automorphism $\nu_\infty$.  Otherwise the proof here proceeds just as in \cite{SelaHopf}.
\end{proof}

Most of the remainder of the paper is devoted to the proof of the following 

\begin{theorem} [cf. Theorem 3.2, p.313, \cite{SelaHopf}] \label{Thm3.2}
All of the edge groups of $\Lambda_{H_\infty^1}$ are finitely generated.
\end{theorem}

Assuming Theorem \ref{Thm3.2}, the argument from pages 313 and 314 of \cite{SelaHopf} proves that $\eta : \Gamma \to \Gamma_\infty$ is an isomorphism, which contradicts what we know, that $\Gamma$ is finitely presented and $\Gamma_\infty$ is not.  This is the contradiction which proves Theorem \ref{GammaHopfian}.  It remains to prove Theorem \ref{Thm3.2}.

\begin{proof}[Proof (Theorem \ref{Thm3.2}).]

Mostly we follow the proof of \cite[Theorem 3.2]{SelaHopf}, and some of the deviations from this proof were inspired by the proof of \cite[Theorem 3.2]{Sela1}.  The strategy is to show that in the presence of edge groups in $\Lambda_{H_\infty^1}$ which are not finitely generated, we can further refine $\Lambda_{H_\infty^1}$, contradicting the canonical properties of the abelian JSJ decomposition.

Note that if not all of the edge groups of $\Lambda_{H_\infty^1}$ are finitely generated then $\eta|_{H_1} : H_1 \to H_\infty^1$ is not an isomorphism, since all abelian subgroups of $H_1$ are finitely generated.  This implies that if there is an edge group of $\Lambda_{H_\infty^1}$ which is not finitely generated then $\psi|_{H_1} : H_1 \to H_1$ is surjective and has nontrivial kernel.

Suppose that $\Lambda_{H_\infty^1}$ does contain edge groups which are not finitely generated.  Erase all the edges with finitely generated stabiliser, and denote by $\Lambda_{L_\infty}$ a remaining connected component which contains an edge.  Let $L_\infty$ be the fundamental group of $\Lambda_{L_\infty}$. 

Let $\gamma_1, \ldots , \gamma_f \in H_1$ be such that $L_\infty = \langle \eta(\gamma_1), \ldots , \eta(\gamma_f) \rangle$.  Once again, by replacing $\psi$ by an iterate composed with an inner automorphism, we may suppose that $L = \langle \gamma_1, \ldots , \gamma_f \rangle$ is invariant under $\psi$ and so $L_\infty$ is the limit group resulting from $\psi : L \to L$.  Since $L$ is a subgroup of $\Gamma$, any abelian subgroup of $L$ is finitely generated.  Hence $\eta |_L : L \to L_\infty$ is not an isomorphism, and $L_\infty$ is a strict $\Gamma$-limit group.  Since $\Gamma$-limits groups are CSA, and $\Lambda_{L_\infty}$ contains more than $1$ vertex, $L_\infty$ is not abelian.

Let $V_\infty^1, \ldots , V_\infty^m$ be the vertex groups and $E_\infty^1, \ldots , E_\infty^s$ the edge groups in $\Lambda_{L_\infty}$. Note that each of the $E_\infty^j$ is not finitely generated.  Let $t_\infty^1, \ldots ,t_\infty^b$ be the set of Bass-Serre generators in the graph of groups $\Lambda_{L_\infty}$.  The vertex groups together with the Bass-Serre generators $t_\infty^j$ generate $L_\infty$. 

Now let $a_1^1, \ldots , a_{l_1}^1, \ldots , a_1^m, \ldots , a_{l_m}^m, t_1, \ldots , t_b$ be elements in $L$ so that  $\eta (a_1^i), \ldots , \eta (a_{l_i}^i) \in V_\infty^i$ for each $1 \le i \le m$, so that $\eta(t_j) = t_\infty^j$ for each $1 \le j \le b$ and so that for each generator $\gamma_j$ of $L$ there is a word:
\[	\eta(\gamma_j) = w_j(\eta(a_1^1), \ldots , \eta(a_{l_m}^m), t_\infty^1, \ldots , t_\infty^b)	.	\]
By possibly replacing the elements $a_1^1, \ldots , a_{l_m}^m , t_1, \ldots , t_b$ by $\psi^k(a_1^1)$, $\ldots$, $\psi^k(a_{l_m}^m), \psi^k(t_1), \ldots , \psi^k(t_b)$ for some $k \ge 1$, we may assume that the elements $a_1^1, \ldots, a_{l_m}^m, t_1, \ldots , t_b$ generate $L$ and that, for each $1 \le i \le m$, the vertex group $V_\infty^i$ is generated by $\eta(a_1^i), \ldots , \eta(a_{l_i}^i)$, together with the edge groups $E_\infty^j$ associated to edges connected to the vertex corresponding to $V_\infty^i$ in $\Lambda_{L_\infty}$.  For each $j$, let $e_0^j, e_1^j, \ldots \in L$ be a set of elements for which $\eta(e_p^j) \in E_\infty^j \smallsetminus \{ 1 \}$ and $E_\infty^j$ is generated by $\{ \eta(e_0^j), \eta (e_1^j) , \ldots \}$.

We now define a sequence $\{ \hat{U}_n \}$ of finitely presented groups together with homomorphisms $\hat{\kappa}_n : \hat{U}_{n-1} \to \hat{U}_n$ which `approximate' the finitely generated group $L$ .  The groups $\hat{U}_n$ admit epimorphisms onto the subgroup $L$, and the decomposition $\Lambda_{L_\infty}$ can be `lifted' to a graph of groups decomposition $\hat\Lambda_n$ of each $\hat{U}_i$ (this decomposition has finitely generated abelian edge groups).  It is this lifting property which allows us to apply the shortening argument in this situation.  The construction of the groups $\hat{U}_n$ is the major innovation in the proof in this paper, otherwise we mostly follow \cite{SelaHopf} and \cite{Sela1}.

To define the groups $\hat{U}_n$ we first follow \cite{SelaHopf} and \cite{Sela1} by defining groups $U_n$, along with epimorphisms $\lambda_n : U_n \twoheadrightarrow L$ and homomorphisms $\kappa_n : U_{n-1} \to U_n$.  The groups $\hat{U}_n$ are then defined in order that we can apply the shortening argument of the previous sections in this case.

We first define the groups $U_n$ iteratively.  Set $U_0$ to be the free group:
\[	U_0 = \langle v_1^1, \ldots , v_{l_1}^1, \ldots , v_1^m, \ldots , v_{l_m}^m , y_1, \ldots , y_b, z_0^1, \ldots , z_0^s \rangle . \]
Define an epimorphism $\lambda_0 : U_0 \to L$ by $\lambda_0(v_p^i) = a_p^i$, $\lambda_0(y_r) = t_r$ and $\lambda_0(z_0^j) = e_0^j$.  Now define the group $U_1$ to be the group generated by elements
\[	\{ v_1^1, \ldots , v_{l_1}^1, \ldots ,v_1^m, \ldots , v_{l_m}^m, y_1, \ldots , y_b, z_0^1, \ldots , z_0^s, z_1^1, \ldots , z_1^s \} ,	\]
subject to the relations which make $\langle z_0^j, z_1^j \rangle \le U_n$ isomorphic to the subgroup generated by $\eta(e_0^j)$ and $\eta(e_1^j)$ in $E_\infty^i$ (under the map which sends $z_r^j$ to $\eta(e_r^j)$, for $r = 1,2$), for each $j = 1 , \ldots , s$.  Because this subgroup of $E_\infty^i$ is a finitely generated abelian group, only finitely many relations are needed.  Clearly there exists a natural homomorphism $\kappa_1 : U_0 \to U_1$.  Since $K_L = \bigcup\limits_{k=1}^\infty \text{ker}(\psi^k)$, for some $d_1 >  0$ there exists an epimorphism $\lambda_1 : U_1 \to L$ defined by $\lambda_1(v_p^i) = \psi^{d_1}(a_p^i)$, $\lambda_1(y_r) = \psi^{d_1}(t_r)$, $\lambda_1(z_0^j) = \psi^{d_1}(e_0^j)$ and $\lambda_1(z_1^j) = \psi^{d_1}(e_1^j)$.  By our definitions of $U_0$, $U_1$, $\kappa_1$, $\lambda_0$ and $\lambda_1$ we have the following commutative diagram:
\[
\begin{CD}
U_0  			@>\kappa_1>> 	U_1\\
@V{\lambda_0}VV				@V{\lambda_1}VV\\
L				@>\psi^{d_1}>>		L
\end{CD}
\]

Having defined $U_0$ and $U_1$ we continue defining the groups $U_n$, the homomorphisms $\kappa_n : U_{n-1} \to U_n$ and the epimorphisms $\lambda_n : U_n \to L$ inductively.  First define the group $G_n$ to be the group generated by
\[	\{	v_1^1, \ldots , v_{l_1}^1, \ldots , v_1^m, \ldots , v_{l_m}^m, y_1, \ldots , y_b, z_0^1, \ldots, z_0^s, \ldots, z_n^1, \ldots , z_n^s	\}	,	\]
together with the (finitely many) relations which force the subgroup $E_n^j = \langle z_0^j, \ldots , z_n^j \rangle$ to be isomorphic, under the map $z_i^j \to \eta(e_i^j)$, to the subgroup of $L_\infty$ generated by $\{ \eta(e_0^j) , \ldots , \eta(e_n^j) \}$, for $j = 1, \ldots , s$.  The group $G_n$ admits a natural epimorphism $\sigma_n : G_n \to L_\infty$, defined by setting $\sigma_n(v_p^i) = \eta(a_p^i)$, $\sigma_n(y_r) = t_\infty^r$ and $\sigma_n(z_d^j) = \eta(e_d^j)$.  The group $U_n$ is a quotient of $G_n$.  To the existing set of relations of $G_n$ we add all words $w$ in the given generating set for $G_n$ for which:
\begin{enumerate}
\item $\sigma_n(w) = 1$;
\item the length of $w$ in the given generating set for $G_n$ is at most $n$; and
\item For some fixed index $i \in \{ 1, \ldots , m \}$, the word $w$ is a word in
\begin{enumerate}
\item The generators $v_1^i, \ldots , v_{l_i}^1$;
\item The elements $z_1^j, \ldots , z_n^j$ for an index $j \in \{ 1, \ldots , s \}$ so that $E_\infty^j < V_\infty^i$; and
\item The words $y_rz_1^jy_r^{-1}, \ldots , y_rz_n^j y_r^{-1}$ for a pair of indices $(j,r)$ with $1 \le j \le s$ and $1 \le r \le b$, for which $t_\infty^rE_\infty^j(t_\infty^r)^{-1} < V_\infty^i$.
\end{enumerate}
\end{enumerate}

Clearly there exists a natural map $\kappa_n : U_{n-1} \to U_n$ and, since $K_L = \bigcup\limits_{k=1}^\infty \text{ker}(\psi^k)$, for some integer $d_n > d_{n-1}$ there exists an epimorphism $\lambda_n : U_n \to L$ defined by $\lambda_n(v_p^i) = \psi^{d_n}(a_p^i)$, $\lambda_n(y_r) = \psi^{d_n}(t_r)$ and $\lambda_n(z_d^j) = \psi^{d_n}(e_d^j)$.  By the definitions of $U_n$, $\kappa_n$ and $\lambda_n$, if we define $k_n = d_n - d_{n-1}$ then the following diagram commutes:
\[
\begin{CD}
U_0		@>\kappa_1>> 	U_1 @. \  \cdots 	@>\kappa_{n-1}>> 	U_{n_1} 	@>\kappa_n>> 	U_n \\
@V\lambda_0VV		@V\lambda_1VV	@.			@V\lambda_{n-1}VV	@V\lambda_nVV \\
L		@>\psi^{d_1}>>		L	@. \  \cdots	@>\psi^{k_{n-1}}>>	L		@>\psi^{k_n}>> L
\end{CD}
\]

Since the second set of defining relations of the group $U_n$ consists of words whose letters are mapped by $\sigma_n$ into the same vertex group in $\Lambda_{L_\infty}$, each of the groups $U_n$ admits an abelian splitting $\Lambda_n$ which projects by $\sigma_n$ into the abelian decomposition $\Lambda_{L_\infty}$ of $L_\infty$.  That is to say each of the vertex groups $V_n^i$ in $\Lambda_n$ satisfies $\sigma_n(V_n^i) \le V_\infty^i$, each of the edge groups $E_n^j = \langle z_1^j, \ldots , z_n^j \rangle$ satisfies $\sigma_n(E_n^j) < E_\infty^j$, and each of the Bass-Serre generators in $\Lambda_n$ satisfies $\sigma_n(y_r) = t_\infty^r$.

We now define the finitely presented groups $\hat{U}_n$.  First define $\hat{U}_0 = U_0$.  We define the group $\hat{U}_n$ using the group $U_n$.  Fix an $n$, and consider the edge groups of $U_n$, namely $E_n^1, \ldots , E_n^s$.  For any $i \in \{ 1, \ldots  , s \}$, since $\eta(e_r^i)$ represents a non-trivial element of $\Gamma_\infty$, we know for all $k \ge 1$ that $\psi^k \circ \lambda_n (E_n^i) \neq \{ 1\}$.  In fact, since $E_n^i$ is isomorphic to $\langle \eta(e_0^i), \ldots , \eta(e_n^i) \rangle \le E_\infty^i$, the homomorphism $\lambda_n$ is injective when restricted to $E_n^i$.  Suppose that $j$ is such that for all $k \ge 0$ there exists an non-cyclic maximal abelian subgroup $M_k$ of $\Gamma$ for which $\psi^{\zeta_k} \circ \lambda_n (E_n^j) \le M_k$, where $\zeta_k = d_{n+k} - d_n$.  

Let $V_n^{p_1}$ and $V_n^{p_2}$ be the vertex groups of $U_n$ adjacent to $E_n^i$.  Suppose that $V_n^{p_1}$ and $V_n^{p_2}$ are both nonabelian and suppose that $\lambda_n(E_n^j) \neq M_0$.  For $s = 1,2$, define $\hat{V}_n^{p_s}$ to be the free product of $V_n^{p_s}$ and $M_0$, amalgamated over the subgroups $E_n^j$ and $\lambda_n(E_n^j)$.  Let $\hat{E}_n^j$ be the copy of $M_0$ in $\hat{V}_n^{p_s}$.

Perform this construction for each pair of nonabelian vertex groups $V_n^{p_1}, V_n^{p_2}$ with such a common edge group $E_n^j$.  If some vertex group $V_n^p$ contains two edge groups of the above kind, $E_n^{j_1}$ and $E_n^{j_2}$ such that $\lambda_n(E_n^{j_1}), \lambda_n(E_n^{j_2}) \le M_0$, where $M_0$ is the same noncyclic maximal abelian subgroup of $\Gamma$, then we embed $E_n^{j_1}$ and $E_n^{j_2}$ in the same subgroup $M_0 \le \hat{V}_n^p$, and let $\hat{E}_n^{j_1} = \hat{E}_n^{j_2} = M_0 \le \hat{V}_n^p$.  We end with a collection $\hat{V}_n^i$ and $\hat{E}_n^j$ of edge and vertex groups (if $E_n^j$ is not an edge group of the above sort, then define $\hat{E}_n^j = E_n^j$).

To define the group $\hat{U}_n$, take the graph of groups decomposition $\Lambda_n$ of $U_n$ described above, and replace each vertex group $V_n^i$ with $\hat{V}_n^i$ and each edge group $E_n^i$ with $\hat{E}_n^j$.  Call the resulting graph of groups $\hat\Lambda_n$.  The group $\hat{U}_n$ is the fundamental group of $\hat\Lambda_n$.

It is not difficult to see that $U_n$ is a subgroup of $\hat{U}_n$, that $\lambda_n$ can be extended in a natural way to $\hat{\lambda}_n : \hat{U}_n \to L$. To define $\hat{\kappa}_n : \hat{U}_{n-1} \to \hat{U}_n$, we need only specify what is done to elements of $U_n$ and what is done to the elements of the groups $M_0$ as above.  Map $U_n$ as prescribed by $\kappa_n$, and map $M_0$ via $\psi^{k_n}$ to $\psi^{k_n} (M_0) \le M_1$, where $M_0$ and $M_1$ are considered as subgroups of $\Gamma$, and $M_1$ is the unique maximal abelian subgroup of $\Gamma$ containing $\psi^{k_n}(M_0)$.  We can now see that the following diagram is commutative:
\[
\begin{CD}
U_0		@>\kappa_1>> 	U_1 @. \  \cdots 	@>\kappa_{n-1}>> 	U_{n_1} 	@>\kappa_n>> 	U_n \\
@V\iota_0VV 			@V\iota_1VV @.	 				@V\iota_{n-1}VV				@V\iota_nVV \\
\hat{U}_0	@>\hat\kappa_1>>\hat{U}_1 @. \ \cdots@>\hat\kappa_{n-1}>>	\hat{U}_{n-1} @>\hat\kappa_n>> \hat{U}_n \\
@V\hat\lambda_0VV	@V\hat\lambda_1VV	@.			@V\hat\lambda_{n-1}VV	@V\hat\lambda_nVV \\
L		@>\psi>>		L	@. \  \cdots	@>\psi^{k_{n-1}}>>	L		@>\psi^{k_n}>> L\ ,
\end{CD}
\]
where the map $\iota_i : U_i \to \hat{U}_i$ denotes inclusion.  Also, $\lambda_i = \hat\lambda_i \circ \iota_i$.

Denote by $\text{Mod}(\hat{U}_n)$ the subgroup of $\text{Aut}(\hat{U}_n)$ generated by inner automorphisms of $\hat{U}_n$, Dehn twists along edges of $\hat\Lambda_n$ and generalised Dehn twists induced by abelian vertex groups in $\hat\Lambda_n$.  

Let $W_n$ be the subgroup of $U_n$ generated by the $v_p^i$ and the $y_r$.  Of course, $W_n$ is also a subgroup of $\hat{U}_n$.  We use the fact that $\lambda_n = \hat\lambda_n \circ \iota_n$ to blur the distinction between these two inclusions of $W_n$.  Clearly, the homomorphism $\kappa_n : U_{n-1} \to U_n$ restricts to an epimorphism from $W_{n-1}$ to $W_n$, and $\lambda_n$ restricts to an epimorphism from $W_n$ onto $L$.  Hence we obtain a sequence of epimorphisms:
\[
\begin{CD}
W_0 @>\kappa_1>> W_1 \ \cdots @>\kappa_{n-1}>> W_{n-1} @>\kappa_n>> W_n @>\kappa_{n+1}>> \ \cdots	\ ,
\end{CD}
\]
where the direct limit of the $\{ W_n \}$ and the maps $\{ \kappa_n \}$ is $L_\infty$.

We now adapt the shortening argument from the first sections of this paper to prove Theorem \ref{Thm3.2}.  Let $\Theta$ be a maximal tree in the graph of groups $\Lambda_{L_\infty}$.  We order the $m$ vertices $p_1, \ldots , p_m$ in $\Theta$ so that $p_1$ is connected to $p_2$, and in general $p_Q$ is connected to the subtree of $\Theta$ spanned by $p_1, \ldots , p_{q-1}$.  Without loss of generality, we assume that this is the original order defined on the vertex groups $V_\infty^i$ of the decomposition $\Lambda_{L_\infty}$.

Let $X$ be the CWIF space upon which $\Gamma$ acts properly, cocompactly and torally by isometries, and let $d_X$ be the metric on $X$.  Let $x \in X$ be the arbitrarily chosen basepoint.  For each $\gamma \in \Gamma$, $n \ge 1$ and $\phi \in \text{Mod}(\hat{U}_n)$, define the stretching constants
\begin{eqnarray*}
\mu_i(n,\gamma,\phi) &  = &  \max_{1 \le j \le l_i} d_X(x, (\tau_\gamma \circ \lambda_n \circ \phi)(v_j^i) . x), \mbox{ and} \\
\chi_r(n,\gamma,\phi) & = & d_X(x,(\tau_\gamma \circ \lambda_n \circ \phi)(y_r) . x)	,
\end{eqnarray*}
where $\tau_\gamma$ is the inner automorphism of $\Gamma$ induced by $\gamma$.  We also define the corresponding $(m+b)$-tuple:
\[	\text{tup}(n,\gamma,\phi) = \left( \mu_1(n,\gamma,\phi), \ldots , \mu_m(n,\gamma,\phi), \chi_1(n,\gamma,\phi), \ldots , \chi_b(n,\gamma,\phi) \right) .	\]
For each $n$, choose $\gamma_n \in \Gamma$ and $\phi_n \in \text{Mod}(\hat{U}_n)$ for which $\text{tup}(n,\gamma_n,\phi)$ is a minimal $(m+b)$-tuple in the set $\{ \text{tup}(n,\gamma , \phi) \}$ with respect to the lexicographic order on $(m+b)$-tuples.

Since the direct limit of the sequence of groups $\{ W_n \}$ and maps $\{ \kappa_n \}$ is $L_\infty$, if $w \in W_0$ is an element for which $\eta \circ \lambda_0 (w) =1$ then there exists some $n_w$ so that for every $n > n_w$ we have $\kappa_n \circ \cdots \circ \kappa_1 (w) = 1$.  let $w_1, w_2 \in W_0$ be a pair of elements for which $\eta \circ \lambda_0(w_1), \eta \circ \lambda_0(w) \in V_\infty^i$ for some $i$.  By the construction of the sequence of groups $\{ \hat{U}_n \}$ and the sequence of splittings $\{ \hat\Lambda_n \}$, for some $n_0$ and every $n > n_0$, the elements $\kappa_n \circ \cdots \circ \kappa_1 (w_1)$ and $\kappa_n \circ \cdots \circ \kappa_1 (w_2)$ belong to the same ($i$th) vertex group $V_n^i$ in the splitting $\hat\Lambda_n$ of $\hat{U}_n$.  Hence, for every $n > n_0$, for every $\phi \in \text{Mod}(\hat{U}_n)$, both $\kappa_n \circ \cdots \circ \kappa_1 (w_1)$ and $\kappa_n \circ \cdots \circ \kappa_1 (w_2)$ are conjugated by the same element in $\Gamma$, and in particular this is true for the chosen $\phi_n$.

In addition, for any index $i$ with $1 \le i \le m$, the subgroup $V_\infty^i$ contains abelian subgroups which are not finitely generated, but there is a (finite) bound on the rank of abelian subgroups of $\Gamma$.  Therefore, there does not exist an infinite sequence of indices, $n_1 < n_2 < \cdots$ and corresponding elements $h_{n_k} \in \Gamma$ so that for every $j$, the element $h_{n_k}$ conjugates $\lambda_{n_j} \circ \phi_{n_j} (v_p^i)$ into $\lambda_{n_1}(v_p^i)$ for $p = 1, \ldots , l_i$.  Thus, the sequence of displacement constants $\{ \| \tau_{\gamma_n} \circ \lambda_n \circ \phi_n \| \}$ (defined with respect to the obvious generating sets of the $W_n$) does not contain a bounded subsequence.

The groups $W_n$ admit a natural isometric action on $X$ via the homomorphism $\tau_{\gamma_n} \circ \lambda_n \circ \phi_n$.  We denote this action by $\rho_n : W_n \to \text{Isom}(X)$.  Since the sequence $\{ \| \rho_n \| \}$ (defined in the obvious way) does not contain a bounded subsequence, and since the group $W_n$ is a natural quotient of the group $W_0$, we may rescale the metric on $X$ by $\frac{1}{\| \rho_n \| }$, and apply the construction from \cite{CWIF1} to find a subsequence (which we still denote by $\{ \rho_n \}$) converging to an action of $W_0$ on a limiting space $\mathcal C_{\infty, \{ \rho_n \}}$, with associated $\R$-tree $T_{\{ \rho_n \}}$.  The actions of $W_0$ on $\mathcal C_{\infty, \{ \rho_n \}}$ and $T_{ \{ \rho_n \} }$ have no global fixed point.

Let $KW_\infty$ be the kernel of the action of $W_0$ on $\mathcal C_{\infty, \{ \rho_n \} }$, and let $Q_\infty = W_0 / KW_\infty$ be the strict $\Gamma$-limit group.  Since $L_\infty$ is non-abelian, so is $Q_\infty$, and we know that $T_{ \{ \rho_n \} }$ is not isometric to a real line.  The conclusions of Theorem \ref{Texists} hold for the actions of $W_0$ and $Q_\infty$ on $T_{\{ \rho_n \} }$.

The group $L_\infty$ is the direct limit of the group $\{ W_n \}$ and the epimorphisms $\{ \kappa_n \}$.  Since if $w_1, w_2 \in W_0$ are a pair of elements for which $\eta \circ \lambda_0 (w_1), \eta \circ \lambda_0(w_2) \in V_\infty^i$ for some $i$, then for some index $n_0$ and every $n > n_0$ the elements $\kappa_n \circ \cdots \circ \kappa_1 (w_1)$ and $\kappa_n \circ \cdots \circ \kappa_1 (w_2)$ belong to the same ($i$th) vertex group $V_n^i$ in the abelian splitting $\Lambda_n$ of $U_n$.  Therefore, for every $n > n_0$, both $\kappa_n \circ \cdots \circ \kappa_1 (w_1)$ and $\kappa_n \circ \cdots \circ \kappa_1 (w_2)$ are conjugated by the same element by the chosen automorphism $\phi_n$.  Hence, $V_\infty^i$ is naturally embedded in the group $Q_\infty$, for each $1 \le i \le m$.  Also, since Dehn twists and generalised Dehn twists fix edge groups elementwise, the intersection of vertex groups $V_\infty^{i_1}, V_\infty^{i_2}$ in $L_\infty$ are embedded as intersections of $V_\infty^{i_1}, V_\infty^{i_2}$ in $Q_\infty$.  Since the maps $\kappa_n : W_{n-1} \to W_n$ and $\lambda_n : W_n \to L$ are epimorphisms for every $n$, the group $Q_\infty$ is generated by the subgroups $V_\infty^i$ and the images of the elements $y_r \in W_0$ in $Q_\infty$.

At this point, because of the possibility of abelian subgroups of $L_\infty$ which are not locally cyclic, we proceed as in \cite[Theorem 3.2]{Sela1}, rather than continuing to follow the proof of \cite[Theorem 3.2]{SelaHopf}.

We first prove that one of the vertex groups $V_\infty^i$ does not fix a point in $T_{\{ \rho_n \}}$.  In order to derive a contradiction, we prove the following.

\begin{proposition} [cf. Proposition 3.3, p.49, \cite{Sela1}] \label{VertexFixBasept}
If all the subgroups $V_\infty^1, \ldots , V_\infty^m$ (the vertex groups in the splitting $\Lambda_{L_\infty}$ of $L_\infty$) fix points in the $\R$-tree $T_{ \{ \rho_n \}}$, then they all fix the basepoint $t_0 \in T_{ \{ \rho_n \} }$.
\end{proposition}
\begin{proof}
Suppose that all of the vertex groups $V_\infty^1, \ldots , V_\infty^m$ fix points in the $\R$-tree $T_{ \{ \rho_n \} }$.  That $V_\infty^1$ fixes the basepoint $t_0 \in T_{ \{ \rho_n \} }$ follows exactly as in the proof of \cite[Proposition 3.6, p. 319]{SelaHopf} or as in the proof of \cite[Proposition 3.3, p.49]{Sela1}.

Suppose then that $V_\infty^1, \ldots , V_\infty^{i_0 -1}$ fix $t_0$ but that $V_\infty^{i_0}$ does not.  Then we can apply the shortening arguments from Section \ref{Discrete} to find $\alpha_n \in \text{Mod}(\hat{U}_n)$ and $\gamma_n' \in \Gamma$ so that for all but finitely many $n$ we have:
\[	\text{tup}(n,\gamma_n',\phi_n \circ \alpha_n) < \text{tup}(n,\gamma_n,\phi_n) .	\]
There are a number of cases to consider.  Note that none of the vertex groups correspond to axial components of $T_{ \{ \rho_n \} }$, since they all fix points in $T_{ \{ \rho_n \} }$.

The two main cases are where the path $\beta$ between $t_0$ and the point fixed by $V_\infty^{i_0}$ is contained in a line $p_E$ and when it is not.  When it is not, we proceed just as in Case 2 from Section \ref{Discrete} above.  When $\beta$ is contained in some $p_E \in T_{ \{ \rho_n \} }$ we use the fact the edge group $E_n^p \in U_n$ corresponding to the edge group $E_\infty^p$ fixing $\beta$ maps {\em onto} the maximal abelian subgroup fixing an approximating flat.  Therefore, we can proceed just as in Section \ref{Discrete}, following \cite[\S 6]{RipsSela1}.
\end{proof}

\begin{lemma} [cf. Lemma 3.5, p. 50, \cite{Sela1}] \label{BassSerreFixBasept}
Suppose that all of the vertex groups $V_\infty^1, \ldots , V_\infty^m$ fix the basepoint $t_0 \in T_{\{ \rho_n \}}$, and that not all of the elements $y_1, \ldots , y_b$ fix $t_0$.  Let $j_0$ be the minimal index for which $y_{j_0}$ does not fix $t_0$.  Then there exists some integer $n_0$ so that for all $n > n_0$ there exist $\beta_n \in \text{Mod}(\hat{U}_n)$ for which:
\begin{enumerate}
\item $\beta_n$ fixes $V_n^i$ elementwise, for $i = 1, \ldots, m$;
\item $\beta_n(y_j) = y_j$ for $j = 1, \ldots , j_0 -1$; and
\item $\chi_{j_0}(n,\gamma_n, \phi_n\circ \beta_n) < \chi_{j_0}(n,\gamma_n,\phi_n)$.
\end{enumerate}
\end{lemma}
\begin{proof}
Again, we proceed as in Section \ref{Discrete}, following \cite[\S 6]{RipsSela1}.
\end{proof}

Proposition \ref{VertexFixBasept} and Lemma \ref{BassSerreFixBasept}, together with the choice of $\phi_n$ and $\gamma_n$, the fact that $Q_\infty$ is generated by the $V_\infty^i$ and the images of $y_1, \ldots , y_b$ and the fact that the action of $Q_\infty$ on $T_{\{ \rho_n \}}$ has no global fixed point imply

\begin{proposition} [cf. Proposition 3.6, p. 51, \cite{Sela1}] \label{NotAllVerticesElliptic}
As least one of the vertex groups in the decomposition $\Lambda_{L_\infty}$ of $L_\infty$ acts non-trivially on the $\R$-tree $T_{\{ \rho_n \}}$.
\end{proposition}

\begin{lemma} \cite[Lemma 3.7, p. 51]{Sela1} \label{AbelianSubgroups}
Let $A$ be an abelian subgroup of one of the groups $V_\infty^i$.  Then either $A$ fixes a point in the $\R$-tree $T_{\{ \rho_n \}}$ or $A$ can be written as a direct sum $A = A_1 \oplus A_2$, where $A_1$ fixes a point in $T_{\{ \rho_n \}}$ and $A_2$ is a finitely generated free abelian group acting freely on $T_{ \{ \rho_n \}}$.
\end{lemma}

\begin{proposition} [cf. Proposition 3.8, p. 51, \cite{Sela1}] \label{AllNonAbVerticesElliptic}
Every nonabelian vertex group $V_\infty^i$ in $\Lambda_{L_\infty}$ fixes a point $t_i \in T_{ \{ \rho_n \}}$.
\end{proposition}
\begin{proof}
Suppose $V_\infty^i$ is nonabelian and does not fix a point in $T_{\{ \rho_n \}}$.  

The argument from pages 317--318 of \cite{SelaHopf} implies that $Q_\infty$ is freely indecomposable rel the edge groups $E_\infty^i$.  Therefore, by Remark \ref{RelativeRtreeStructure}, the action of $Q_\infty$ on $T_{ \{ \rho_n \} }$ induces a graph of groups decomposition of $Q_\infty$ with abelian edge groups.  Using Lemma \ref{AbelianElliptic} we can modify this graph of groups to get a non-trivial abelian splitting $\Omega$ of $Q_\infty$ in which all non-cyclic maximal abelian subgroups are elliptic.  Since $V_\infty^i$ is a subgroup of $Q_\infty$, the graph of groups $\Omega$ induces a graph of groups decomposition $\Omega^i$ of $V_\infty^i$ with all non-cyclic abelian subgroups elliptic.  Since we are assuming that $V_\infty^i$ does not fix a point in $T_{ \{ \rho_n \} }$, the graph of groups decomposition $\Omega^i$ is non-trivial (applying Lemma \ref{AbelianElliptic} does not change the vertex groups).

Let $T_{\Omega^i}$ be the Bass-Serre tree of the splitting $\Omega^i$ of $V_\infty^i$.  Consider the edge groups in the abelian JSJ decomposition of $H_\infty^1$ which are adjacent to $V_\infty^i$.  By the above, any non-cyclic edge group is elliptic in this splitting.  However, the conjugacy class of a generator of any cyclic edge group is periodic under $\nu_\infty$, so is also elliptic in $T_{\Omega^i}$.  Therefore, the action of $V_\infty^i$ on $T_{ \Omega^i}$ can be extended to a (stable) action with trivial stabilisers of tripods of $H_\infty^i$ on some $\R$-tree $\hat{T}$ by letting all the other vertex groups in the abelian JSJ decomposition of $H_\infty^1$ fix points in $\hat{T}$.  By Theorem \ref{RtreeStructure}, $H_\infty^1$ inherits an abelian splitting $\hat\Lambda$ from its action on the $\R$-tree $\hat{T}$.  Since $V_\infty^i$ does not fix a point in $T_{\Omega^i}$, it does not fix a point in $\hat{T}$, so $V_\infty^i$ can not be conjugated into a vertex group in $\hat\Lambda$.  The group $V_\infty^i$ contains abelian subgroups which are not finitely generated, so is not a CMQ subgroup of $H_\infty^i$.  However, this contradicts the properties of the JSJ decomposition (see Theorem \ref{AbelianJSJ}.(\ref{JSJ-Vertex})).  This contradiction proves the proposition.
\end{proof}

By Propositions \ref{NotAllVerticesElliptic} and \ref{AllNonAbVerticesElliptic}, there must be some abelian vertex group $V_\infty^{i_0}$ which does not fix a point in $T_{ \{ \rho_n \}}$.  By Lemma \ref{AbelianSubgroups}, $V_\infty^{i_0} = A^{i_0} \oplus D^{i_0}$, where $A^{i_0}$ fixes a point in the $\R$-tree $T_{\{ \rho_n \}}$ and $D^{i_0}$ is a non-trivial, finitely generated free abelian group that acts faithfully on the axis on $V_\infty^{i_0}$.  Since every maximal abelian subgroup in $L_\infty$ is malnormal, all adjacent vertex groups to an abelian vertex group are non-abelian.  Since all non-abelian vertex groups fix points in $T_{ \{ \rho_n \}}$, all edge groups connected to $V_\infty^{i_0}$ are subgroups of $A^{i_0}$.  Note that since $V_\infty^i$ contains an infintely generated abelian subgroup (by the construction of $L_\infty$), the subgroup $A^{i_0}$ is itself infinitely generated.

Let $L_\infty^1$ be the fundamental group of the graph of groups $\Upsilon_{L_\infty^1}$ obtained from $\Lambda_{L_\infty}$ by replacing the vertex group $V_\infty^{i_0}$ with its subgroup $A^{i_0}$.  Clearly $L_\infty^1$ is a subgroup of $L_\infty$, and since $V_\infty^{i_0} = A^{i_0} \oplus D^{i_0}$, the group $L_\infty^1$ is also a quotient of $L_\infty$.  Hence, $L_\infty^1$ is finitely generated.  Also, the first Betti number of $L_\infty^1$ is strictly smaller than that of $L_\infty$: $b_1(L_\infty^1) < b_1(L_\infty)$. 

We can repeat the above shortening arguments for the group $L_\infty^1$, along with its graph of groups decomposition $\Upsilon_{L_\infty^1}$.  This gives us a new $\Gamma$-limit group $Q_\infty^1$ acting on a new $\R$-tree $T^1$.  By the above arguments, one of the abelian vertex groups of $\Upsilon_{L_\infty^1}$ acts nontrivially on $T^1$, and we may obtain a new $\Gamma$-limit group $L_\infty^2$ by the same procedure as that which found $L_\infty^1$ from $L_\infty$.  Also $b_1(L_\infty^2) < b_1(L_\infty^1) < b_1(L_\infty)$.  Therefore, after repeating this argument finitely many times, we finally arrive at a contradiction, which concludes the proof of Theorem \ref{Thm3.2} and hence also that of Theorem \ref{GammaHopfian}.
\end{proof}


\begin{thebibliography}{99}
\bibitem{BF2} M. Bestvina and M. Feighn, Bounding the complexity of simplicial group actions, \textit{Invent. Math.} {\bf 103} (1991), 449--469.
\bibitem{BF} M. Bestvina and M. Feighn, Stable actions of groups on real trees, \textit{Invent. Math.} {\bf 121} (1995), 287--321.
\bibitem{BFSela} M. Bestvina and M. Feighn, Notes on Sela's work: Limit groups and Makanin-Razborov diagrams, preprint.
\bibitem{Dahmani} F. Dahmani, Accidental parabolics and relatively hyperbolic groups, preprint.
\bibitem{DP} T. Delzant and L. Potyagailo, Endomorphisms of Kleinian groups, \textit{GAFA} {\bf 13} (2003), 396--436.
\bibitem{DS} C. Dru\c{t}u and M. Sapir, Tree-graded spaces and asymptotic cones of groups, preprint.
\bibitem{DunSag} M. Dunwoody and M. Sageev, JSJ-splittings for finitely presented groups over slender groups, \textit{Invent. Math.} {\bf 135} (1999), 25--44.
\bibitem{FP} K. Fujiwara and P. Papasoglu, JSJ decompositions of finitely presented groups and complexes of groups, preprint.
\bibitem{CWIF1} D. Groves, Limits of (certain) CAT$(0)$ groups, I: Compactification, preprint.
\bibitem{RelHypCSA} D. Groves, Limit groups for relatively hyperbolic groups, in preparation.
\bibitem{coHopf} D. Groves, On the co-Hopf property for relatively hyperbolic groups, in preparation.
\bibitem{Hruska} C. Hruska, Nonpositively curved spaces with isolated flats, PhD thesis, Cornell University, 2002.
\bibitem{HruskaEmail} C. Hruska, private communication, 2004.
\bibitem{KL} M. Kapovich and B. Leeb, On asymptotic cones and quasi-isometry classes of fundamental groups of $3$-manifolds, \textit{GAFA} {\bf 5} (1995), 582--603.
\bibitem{Levitt} G. Levitt, La dynamique des pseudogroupes de rotations, \textit{Invent. Math.} {\bf 113} (1993), 633--670.
\bibitem{Morgan} J. Morgan, Ergodic theory and free actions of groups on $\R$-trees, \textit{Invent. Math.} {\bf 94} (1988), 605--622.
\bibitem{OP} K. Oshika and L. Potyagailo, Self-embeddings of Kleinian groups, \textit{Ann. Ec. Norm. Sup.} {\bf 31} (1998), 329-343.
\bibitem{Paulin3} F. Paulin, Topologie de Gromov \'equivariante, structures hyperboliques et arbres r\'eels, {\textit Invent. Math.} {\bf 94} (1988), 53--80.
\bibitem{RipsSela1} E. Rips and Z. Sela, Structure and rigidity in hyperbolic groups, I, \textit{GAFA} {\bf 4} (1994), 337--372.
\bibitem{RipsSela} E. Rips and Z. Sela, Cyclic splittings of finitely presented groups and the canonical JSJ decomposition, \textit{Ann. Math.} {\bf 146} (1997), 53--104.
\bibitem{SS} P. Scott and G.A. Swarup, Regular neighbourhoods and canonical decompositions of groups, \textit{Ast\'erisque}, {\bf 289} (2003).
\bibitem{SelaGAFA} Z. Sela, Structure and rigidity in (Gromov) hyperbolic groups and discrete groups in rank $1$ Lie groups II, \textit{GAFA} {\bf 7} (1997), 561--593.
\bibitem{SelaAcyl} Z. Sela, Acylindrical accessibility for groups, \textit{Invent. Math.} {\bf 129} (1997), 527--565.
\bibitem{SelaHopf} Z. Sela, Endomorphisms of hyperbolic groups I: The Hopf property, \textit{Topology} {\bf 38} (1999), 301--321.
\bibitem{Sela1} Z. Sela, Diophantine geometry over groups, I. Makanin-Razborov diagrams. \textit{Pbul. Math. IHES} {\bf 93} (2003),  31--105. 
\bibitem{SelaProblems} Z. Sela, Diophantine geometry over groups: a list of research problems, \texttt{http://www.ma.huji.ac.il/$\sim$zlil/problems.dvi/}
\bibitem{Weidmann} R. Weidmann, The Nielsen method for groups acting on trees, \textit{Proc. London Math. Soc. (3)} {\bf 85} (2002), 93--118.
\bibitem{Wise} D. Wise, A non-Hopfian automatic group, \textit{J. Alg.} {\bf 180} (1996), 845--847.
\end{thebibliography}
\end{document}